\magnification=1200
\documentstyle{amsppt}
\hoffset=-0.5pc
\nologo
\vsize=57.2truepc
\hsize=38.5truepc

\spaceskip=.5em plus.25em minus.20em

\define\fra{\frak}
\define\Bobb{\Bbb}
\define\almolone{1}
\define\atiyaone{2}
\define\cartanon{3}
\define\cartanse{4}
\define\carteile{5}
\define\cheveile{6}
\define\duponboo{7}
\define\eilmacon{8}
\define\grehalva{9}
 \define\herzone{10}
\define\poiscoho{11}
\define\souriau{12}
\define\exteweil{13}
\define\perturba{14}
\define\husmosta{15}
\define\kamtonth{16}
\define\kobanomi{17}
\define\koszulon{18}
\define\koszultw{19}
\define\maclaboo{20}
\define\mackbook{21}
 \define\macktwo{22}
 \define\mackthr{23}
 \define\maybook{24}
\define\milnstas{25}
\define\munkholm{26}
  \define\osborn{27}
\define\palaione{28}
\define\pradione{29}
\define\praditwo{30}
\define\pradithr{31}
\define\pradifou{32}
\define\rinehart{33}
\define\stashtwo{34}
\define\stashfiv{35}
\define\stashnin{36}
 \define\teleman{37}
\define\whitnone{38}
\define\whitntwo{39}

\define\Ext{\text{\rm Ext}}
\define\MEExt#1#2#3#4{\Ext^{#1}_{#2}({#3},{#4})}

\noindent
dg-ga/9706002

\topmatter
\title Extensions of Lie-Rinehart algebras and the
Chern-Weil construction
\endtitle
\author Johannes Huebschmann\endauthor
\affil 
Universit\'e des Sciences et Technologies
de Lille
\\
U. F. R. de Math\'ematiques
\\
F-59 655 VILLENEUVE D'ASCQ, France
\\
Johannes.Huebschmann\@univ-lille1.fr
\endaffil
\date{May 17, 1997}
\enddate
\address{
Universit\'e des Sciences et Technologies
de Lille
\newline
U. F. R. de Math\'ematiques
\newline
F-59 655 VILLENEUVE D'ASCQ, France
\newline
Johannes.Huebschmann\@univ-lille1.fr
}
\endaddress
\endkeywords
\endsubjclass
\abstract{
A Chern-Weil construction for
extensions
of Lie-Rinehart algebras is introduced.
This generalizes the classical Chern-Weil construction
in differential geometry
and yields 
characteristic classes for arbitrary extensions
of Lie-Rinehart algebras.
Some examples arising
from spaces with singularities and from foliations
are given that cannot be treated 
by means of the classical Chern-Weil construction.}
\endabstract
\endtopmatter
\document
\leftheadtext{Johannes Huebschmann}
\rightheadtext{Chern-Weil construction}

\beginsection Introduction

Let $R$ be a commutative ring, and let $A$ be a commutative $R$-algebra.
We introduce a general Chern-Weil construction
which yields characteristic classes for extensions of
$(R,A)$-Lie algebras.
These $(R,A)$-Lie algebras have been introduced by
{\smc Herz}~\cite\herzone\ 
under the name 
\lq\lq pseudo-alg\`ebre de Lie\rq\rq\ 
and were examined by {\smc Palais}~\cite\palaione\ 
under the name \lq\lq $d$-Lie ring\rq\rq\ 
and thereafter by {\smc Rinehart}~\cite\rinehart,
who introduced the terminology
\lq\lq $(R,A)$-Lie algebra\rq\rq.
An $(R,A)$-Lie algebra is a Lie algebra $L$ over the ground ring $R$,
together with 
an action of $L$ on $A$ and
an $A$-module structure on $L$, and 
the two structures
satisfy suitable compatibility conditions which generalize
the usual properties of the Lie algebra of smooth vector fields
on a smooth manifold viewed as a module over its ring of smooth
functions; a precise definition will be reproduced in Section 1 below.
With a suitable notion of morphism, such pairs $(A,L)$
constitute a category,
and 
we shall refer to such a pair as a
{\it Lie-Rinehart algebra\/}.
\smallskip
A principal bundle
gives rise to an
extension of
Lie-Rinehart algebras
which arise as spaces of sections of an extension
of vector bundles,
introduced
by {\smc Atiyah}~\cite\atiyaone\
and now usually called
the {\it Atiyah sequence\/} of the principal bundle;
see (2.2) below for details.
It is also common to talk about
{\it transitive Lie algebroids\/},
cf. \cite{\almolone,\ \mackbook,
\pradione--\pradifou}.
{\smc Pradines\/} \cite{\pradione,\ \praditwo}\ 
in fact
introduced the more general concept of (not necessarily
transitive) Lie algebroid,
but the general notion does not involve
an extension in the sense studied in this paper.
{\smc Almeida and Molino} 
\cite\almolone\  
have shown
that Lie's third theorem does not hold
for transitive Lie algebroids:
{\it not every
transitive Lie algebroid
integrates to a principal bundle\/}.
This provides a negative answer to a question
raised by Pradines. 
{\smc Mackenzie}~\cite\mackbook\ developed obstructions
for the integrability of
Lie algebroids.
\smallskip
For a principal bundle,
our Chern-Weil construction boils down
to the usual Chern-Weil construction
~\cite\cartanon,
~\cite\duponboo,
~\cite\kobanomi,
~\cite\osborn;
details will be given in Section 3 below.
However, there are interesting examples, e.~g. 
arising from foliations or from quantization problems, that do
not come from a principal bundle; we shall describe some
such examples in Section 4 below.
A Chern-Weil homomorphism for a transitive Lie algebroid
has been set up by {\smc Teleman}~\cite\teleman,
and our construction extends that of Teleman as well.
The \lq\lq transverse Chern-Weil map\rq\rq\ for
equivariant principal bundles
constructed in \cite\macktwo\ 
provides a Chern-Weil theory for extensions of principal
bundles.
\smallskip
We now give a brief overview of the contents of the paper:
In Section 1  we  
recall some of the basic notions.
In Section 2
we generalize the 
usual concepts of
connection
and curvature
in a principal bundle
to arbitrary
extensions
of Lie-Rinehart algebras.
Within the
category of
finite rank smooth vector bundles,
these notions have already been rephrased
by {\smc Mackenzie}~\cite\mackbook\ 
in the language
of 
Atiyah sequences and
transitive Lie algebroids.
In Section 3 we
introduce a Chern-Weil construction for an extension 
$$
0
@>>>
L'
@>>>
L
@>>>
L''
@>>>
0
\tag0.1
$$
of Lie-Rinehart algebras under the assumption that
the extension splits in the category of $A$-modules.
This Chern-Weil construction furnishes
a morphism
$$
\roman{Hom}_A(\Sigma'_A[s^2L'],A)^L
\longrightarrow
\roman{Alt}_A(L'',A)
\tag0.2
$$
of
differential graded commutative $R$-algebras
whose induced morphism
on homology depends only on the 
congruence class of the
extension
{\rm (0.1)};
see (3.8.1) and (3.8.2) below for details.
Here 
$s^2L'$
refers to the double suspension
of $L'$, and $\Sigma'_A[s^2L']$
is the {\it symmetric coalgebra\/} on $s^2L'$ over $A$,
so that its dual
$\roman{Hom}_A(\Sigma'_A[s^2L'],A)$
is a graded commutative algebra;
furthermore,
as usual 
the notation
${-}^L$ indicates the invariants with respect to the induced
$L$-action, and
the source 
$\roman{Hom}_A(\Sigma'_A[s^2L'],A)^L$
is equipped with the zero differential.
The morphism (0.2) involves the notion of 
curvature for an extension of the kind (0.1),
and in order for this curvature to be defined,
the underlying extension of $A$-modules must split.
When $L'$ is finitely generated and free as an $A$-module,
the graded algebra 
$\roman{Hom}_A(\Sigma'_A[s^2L'],A)$ is just the polynomial $A$-algebra
on an $A$-basis of the dual of $s^2L'$
but the above description
works {\it without any finiteness assumption \/}.
\smallskip
In view of its complete generality,
our approach is likely to
prove useful
for 
geometrical systems 
in infinite dimensions and, furthermore,
for 
systems 
with singularities where 
e.~g. \lq\lq smooth\rq\rq\ functions are to be understood in the sense of 
{\smc Whitney\/}~\cite\whitnone,~\cite\whitntwo,
see Section 4 for details.
While
{\it in the classical finite dimensional
case 
without singularities
there was no need to distinguish\/} between
e.~g. formal differentials and differential forms
and hence not between an 
exterior algebra $\Lambda_A \fra g^*$
over the dual of a Lie algebra $\fra g$ over an algebra $A$
and an algebra
$\roman{Alt}_A(\fra g,A)$
of differential forms and likewise, between
a symmetric algebra
$\Sigma_A \fra g^*$
and a corresponding algebra
$\roman{Hom}_A(\Sigma_A'[\fra g],A)$
of forms etc.,
under more general circumstances
when the requisite $A$-modules are no longer projective
or not even reflexive
more care is needed
and the appropriate objects to work with are those
involving $A$-valued forms etc.
This claim is well illustrated
by our Chern-Weil map
where the {\it true\/} algebra of characteristic
classes
for an extension (0.1) of 
Lie-Rinehart algebras
is
an algebra of the kind
$\roman{Hom}_A(\Sigma'_A[s^2L'],A)^L$
rather than the symmetric algebra
on the dual of $s^2L'$.
See in particular what is said at the beginning of Section 3.
\smallskip
I am much indebted to K. Mackenzie for a number of most valuable 
comments
on a draft of the paper.
It is a pleasure to dedicate this paper to Jim Stasheff;
it has in fact been influenced by
his work on 
characteristic classes
\cite\milnstas\ and on
differential homological algebra,
cf. e.~g. \cite\husmosta.
Furthermore,
about ten years ago, 
studying
his paper \cite\stashtwo,
I 
encountered
Lie-Rinehart algebras for the first time
and thereafter discovered their significance
for Poisson structures \cite{\poiscoho,~\souriau}. 
Since then Jim encouraged me to 
study these notions
per se,
followed my investigations
and generously offered support.
Lie-Rinehart algebras and strong homotopy generalizations thereof
also arose in his work on homological perturbations
in the BRST-description of constrained
hamiltonian systems and variants thereof;
for these matters, see e.~g. his papers \cite{\stashtwo--\stashnin}. 
In this area, there is  still a huge unexplored territory.

\beginsection 1. Lie-Rinehart algebras and their modules

We review briefly the concept of
Lie-Rinehart algebras.
We then  give descriptions of 
appropriate module-,
algebra-, coalgebra-, Lie algebra-, etc.
structures over Lie-Rinehart algebras
and spell out the corresponding differential graded
objects.
\smallskip
Let $R$ be a commutative ring,
fixed throughout;
the unadorned tensor product symbol
$\otimes$ will
always refer to the tensor product over $R$.
Further, let 
$A$ be a commutative
$R$-algebra.
An
$(R,A)$-{\it Lie algebra\/}
\cite\rinehart\ 
is a Lie algebra $L$ over $R$
which acts on (the left of $A$) by derivations
(written ${(\alpha \otimes a) \mapsto \alpha a}$),
and is also an $A$-module
(the structure map being written ${(a\otimes \alpha) \mapsto a\alpha}$),
 in such a way that
suitable compatibility conditions are satisfied
which generalize the usual properties
of the Lie algebra of vector fields
on a smooth manifold
viewed as a module over its ring of functions;
these conditions read
$$
\alignat 2
(a\,\alpha)(b) &= a\,(\alpha(b)),
\quad &&\alpha \in L, \,a,b \in A,
\tag1.1.a
\\
\lbrack 
\alpha,
a\,\beta
\rbrack
&=
a\,
\lbrack 
\alpha,
\beta
\rbrack
+
\alpha(a)\,\beta,\quad &&\alpha,\beta \in L, \,a \in A.
\tag1.1.b
\endalignat
$$
When the emphasis is on the pair
$(A,L)$,
with the mutual structure of interaction,
we refer to a {\it Lie-Rinehart algebra\/}.
Given 
two 
Lie-Rinehart algebras
$(A,L)$
and
$(A',L')$,
a {\it morphism \/}
${
(\phi,\psi)
\colon
(A,L)
\longrightarrow
(A',L')
}$
of 
{\it Lie-Rinehart algebras\/}
is 
the obvious thing, that is,
$\phi$ and $\psi$ are morphisms in the appropriate categories
that are compatible with the additional structure.
With this notion of
morphism, 
Lie-Rinehart algebras
constitute a category.
Apart from the example of smooth functions and smooth
vector fields on a smooth manifold,
a related (but more general) example 
is the pair consisting of a commutative algebra
$A$
and the $R$-module $\roman{Der}(A)$ of derivations of 
$A$
with the obvious $A$-module structure;
here the commutativity of $A$ is crucial.
\smallskip
Let $L$ be an
$(R,A)$-Lie algebra. An
$R$-module $M$
having the structures of a
left $A$-module and 
that of 
a left $L$-module
$\omega
\colon
L
\to \roman{End}(M)
$
is said to be an $(A,L)$-{\it module\/} provided
the actions are compatible,
that is to say, for 
$\alpha
 \in L,\,a \in A,\,m \in M,$
$$
\align
(a\,\alpha)(m) &= a(\alpha(m)),
\tag1.2.a
\\
\alpha(a\,m) &=   a\,\alpha(m) + \alpha(a)\,m.
\tag1.2.b
\endalign
$$
For example,
given a smooth manifold,
such a structure 
with respect to
the Lie algebra of smooth vector fields
on the module of smooth sections of
a smooth vector bundle
is just  a flat connection.
\smallskip
Let $L$ be an
$(R,A)$-Lie algebra and
$M$ an $(A,L)$-module.
The
$R$-multilinear
alternating functions from $L$ into $M$
with the {\smc Cartan-Chevalley-Eilenberg\/}~\cite\cheveile\ 
differential $d$ given by 
$$
\aligned
(df)(\alpha_1,\dots,\alpha_n)
&=
\quad 
(-1)^n
\sum_{i=1}^n (-1)^{(i-1)}
\alpha_i(f (\alpha_1, \dots\widehat{\alpha_i}\dots, \alpha_n))
\\
&\phantom{=}+\quad
(-1)^n
\sum_{j<k} (-1)^{(j+k)}f(\lbrack \alpha_j,\alpha_k \rbrack,
\alpha_1, \dots\widehat{\alpha_j}\dots\widehat{\alpha_k}\dots,\alpha_n)
\endaligned
\tag1.3
$$
constitute a chain
complex
$\roman{Alt}_R(L,M)$
where
as usual \lq $\ \widehat {}\ $\rq \ 
indicates omission of the corresponding term.
The sign $(-1)^n$ in (1.3) has been introduced 
according to 
the usual Eilenberg-Koszul
convention in differential homological algebra
for
consistency with 
(1.3$'$) below
and with 
what is said in our follow up paper~\cite\exteweil;
in the classical approach such a sign does not occur.
As observed first by
{\smc Palais}~\cite\palaione,
the differential $d$ on
$\roman{Alt}_R(L,M)$
passes to an
$R$-linear
differential on the
graded $A$-submodule
$\roman{Alt}_A(L,M)$
of $A$-multilinear functions,
written $d$, too;
this differential will not be
 $A$-linear unless $L$ acts trivially on $A$, though.
Before we proceed further we mention that
a {\it distinction\/} between {\it graded\/} $A$-objects
and {\it differential graded\/} $R$-objects will persist throughout.
We shall carry out most constructions, e.~g. coalgebras, algebras, etc.
over $A$; however, in view of the non-triviality of the action of $L$ 
on $A$,
most resulting differential graded objects will
be over the ground ring $R$ only.
\smallskip
Given $(A,L)$-modules
$M'$ and $M''$,
the usual formula 
$$
\alpha(x \otimes_A y) = \alpha(x) \otimes_A y + x \otimes_A \alpha(y),
\quad
\alpha \in L, 
\,x \in M',\, y \in M'',
\tag1.4.1
$$
endows
the tensor product $M' \otimes_A M''$ with the 
structure of an $(A,L)$-module,
referred to
as the
{\it tensor product of\/}
$M'$ and $M''$
{\it in the category of
$(A,L)$-modules\/};
if
$M$
is another $L$-module, 
a pairing 
$
\mu_A
\colon
M'
\otimes_A
M''
\longrightarrow
M
$
of $A$-modules
which is a morphism of $(A,L)$-modules
(with respect to (1.4.1)) will be said to be a
{\it a pairing of $(A,L)$-modules\/}.
Given such a pairing
$\mu_A$  of $(A,L)$-modules,
let
$${
\mu = \mu_A \pi
\colon
M'
\otimes_R
M''
\longrightarrow
M'
\otimes_A
M''
\longrightarrow
M
}$$
be the indicated pairing of
$L$-modules, and
{\it define\/}
the  shuffle multiplication of $R$-multilinear,
alternating maps 
by
$$
\aligned
(\alpha \wedge \beta)&(x_1,\dots,x_{p+q}) 
\\
&= (-1)^{|\alpha||\beta|}
\sum_{\sigma}\roman{sign}(\sigma)
\mu(\alpha(x_{\sigma(1)},\dots,x_{\sigma(p)})
\otimes \beta(x_{\sigma(p+1)},\dots,x_{\sigma(p+q)})),
\endaligned
\tag1.4.2
$$
where $\alpha \in 
\roman{Alt}^p_R(L,M'), \ \beta \in \roman{Alt}^q_R(L,M''),
\ x_1,\dots,x_{p+q} \in L$. This
yields a pairing
$$
\wedge
\colon
\roman{Alt}_R(L,M') \otimes
\roman{Alt}_R(L,M'')
\longrightarrow
\roman{Alt}_R(L,M)
\tag1.5
$$
of chain complexes
 which is associative in the obvious sense;
here $\sigma$
runs through $(p,q)$-shuffles
and $\roman{sign}(\sigma)$ refers to the sign of $\sigma$.
The sign $(-1)^{|\alpha||\beta|}$
in the formula (1.4.2)
does usually {\it not\/} occur in the descriptions given
in the literature.
This 
sign is dictated
by the graded tensor product in
differential homological algebra, and since
we shall have occasion to extend
the pairing (1.5) to the differential graded context
we must insist on this sign.
See (1.6.5) below for details.
The pairing (1.5)
induces a pairing
$$
\wedge\colon
\roman{Alt}_A(L,M') \otimes_R
\roman{Alt}_A(L,M'')
\longrightarrow
\roman{Alt}_A(L,M)
\tag1.5$'$
$$
of chain complexes over $R$,
still denoted by
$\wedge$.
In particular,
the graded commutative $A$-algebra
${\roman{Alt}_A(L,A)}$
inherits a structure of a differential graded commutative algebra
over the ground ring $R$
but {\it not\/} over $A$ unless $L$ acts trivially on $A$.
The pairing
$\wedge$
plainly factors through a pairing
$$
\roman{Alt}_A(L,M') \otimes_A
\roman{Alt}_A(L,M'')
\longrightarrow
\roman{Alt}_A(L,M)
\tag1.5$''$
$$
of graded $A$-modules,
uniquely determined by the given data,
but (1.5'') will {\it not\/}
be compatible with the differentials
unless $L$ acts trivially on $A$.
Another description of the pairing (1.5'') will be given in (1.6.5) below.
\smallskip
A conceptual explanation 
of  these facts 
in terms of standard homological algebra
over a suitable universal algebra
is 
due to 
{\smc Rinehart}~\cite\rinehart\ 
and has been elaborated upon in our paper~\cite\poiscoho\ 
to which we refer for details.
Here we only mention that
any Lie-Rinehart algebra $(A,L)$ determines
a universal $R$-algebra $U(A,L)$; for example,
when $A$ is the algebra of smooth functions on a smooth
manifold $N$ and  $L$  the Lie algebra of 
smooth vector fields
on $N$, then $U(A,L)$ is the {\it algebra of\/}
(globally defined)
{\it differential operators
on\/} $N$.
For  an $(A,L)$-module $M$,
the {\it cohomology
of $L$ with coefficients in\/} $M$
is then defined by
$$
\roman H_A^*(L,M) = \MEExt*{U(A,L)}AM,
$$
cf. \cite{\poiscoho,\ \rinehart}.
When $L$ is projective as an $A$-module,
the chain complex
$(\roman{Alt}_A(L,M),d)$
(reproduced above)
computes this cohomology;
in the present paper we
shall exclusively
work with the chain complex
$(\roman{Alt}_A(L,M),d)$,
whether or not $L$ is
projective as an $A$-module.
\medskip
\noindent
{\smc 1.6. More structure.}
\smallskip
\noindent
Let $L$ be an $(R,A)$-Lie algebra,
and let
$M$ be an $(A,L)$-module,
with structure map
$\omega \colon L \to \roman{End}(M)$.
At times we shall 
assume $M$ equipped with additional structure,
e.~g. that of an algebra, a chain complex,
a Lie algebra, etc.
Chain complexes with a non-zero differential
will not explicitly  occur
as $(A,L)$-modules, though; all we shall 
need are graded $(A,L)$-modules 
but we can handle chain complexes at no extra cost and hence
we shall do so.
The material to be given until the end
of this Section  is mostly
folk-lore;
since it it is difficult to
give precise references,
we explain some of the requisite details.
We note that
$A$ could be just the ground ring $R$, but in view of later applications
it will be convenient to distinguish
carefully between $A$ and $R$.
\smallskip
\noindent
{\smc 1.6.1. Lie algebras over $A$.}
Let $M=\fra g$ be a Lie algebra over $A$,
with Lie bracket $[\cdot,\cdot] \colon \fra g \otimes_A \fra g \to \fra g$.
We shall 
say that $\fra g$ is
an $(A,L)$-{\it Lie algebra\/}
if 
the structure map
$[\cdot,\cdot]$
is a morphism of $(A,L)$-modules
or, what amounts to the same,
if the values of
$\omega$ lie in
$\roman{Der}(\fra g) \subseteq \roman{End}(\fra g)$;
we shall then occasionally write
$\omega \colon L \to \roman{Der}(\fra g)$.
For an
$(A,L)$-Lie algebra $\fra g$,
the pairing
(1.5$'$) with respect to
 $\mu_A=[\cdot,\cdot]$
endows
the $A$-multilinear forms $\roman{Alt}_A(L,\fra g)$
with values in $\fra g$
with the structure of a 
differential
graded
Lie algebra
over $R$.
In the special case where $A=R$ with trivial action this is of course well 
known.
\smallskip
\noindent
{\smc 1.6.2. 
Chain complexes over $A$.}
Let $M=C$ be a chain complex over $A$,
and write $Z^0(\roman{End}(C))$
for the (homogeneous) chain maps from $C$ to itself
(of degree zero),
i.~e. for the cycles in the corresponding Hom-complex.
We shall say that
$C$ is an $(A,L)$-{\it chain complex\/} provided the action
$\omega$
factors through
an action
$\omega \colon L \to Z^0(\roman{End}(C))$
(denoted still by $\omega$ with a slight abuse of notation).
For an
$(A,L)$-chain complex $C$,
the 
differential
(1.3)
endows
the $A$-multilinear forms $\roman{Alt}_A(L,C)$
with values in $C$
with the structure of a 
chain complex
over $R$,
but some more care is needed to explain what this really means:
Write 
$$
d^0 \colon
\roman{Alt}_A(L,C)
\longrightarrow
\roman{Alt}_A(L,C)
$$
for the differential induced by the differential on $C$,
and write
$$
d^1 \colon
\roman{Alt}_A(L,C)
\longrightarrow
\roman{Alt}_A(L,C)
$$
for the corresponding operator
(1.3), interpreted suitably with respect to the grading;
in other words, 
for a {\it homogeneous\/}
multilinear
alternating
form $f$
on $L$ with values in $C$
of degree $(n-1)$,
i.~e. a sequence
$f =\{f_{n-1},f_{n}, \dots, f_{n+\ell} ,\dots\}$
of 
multilinear,
alternating
forms $f_{n+\ell}$
in $n+\ell$
variables
on $L$ 
with values in $C_{\ell +1}$,
where $\ell \geq -1$,
define
$d^1$ by
$$
\aligned
(d^1f)(\alpha_1,\dots\alpha_n)
&=
\quad 
(-1)^n
\sum_{i=1}^n (-1)^{(i-1)}
\alpha_i(f (\alpha_1, \dots\widehat{\alpha_i}\dots, \alpha_n))
\\
&\phantom{=}+\quad
(-1)^n
\sum_{j<k} (-1)^{(j+k)}f(\lbrack \alpha_j,\alpha_k \rbrack,
\alpha_1, \dots\widehat{\alpha_j}\dots\widehat{\alpha_k}\dots,\alpha_n).
\endaligned
\tag1.3$'$
$$
Then
$d^1$ endows
$\roman{Alt}_A(L,C)$
with the structure of a chain complex,
and it is a standard fact that so does $d^0$.
Moreover,
since the 
action of $L$ on $C$
is assumed compatible with the differential
on $C$,
the operator $d^1$
is also compatible
with the differential
$d^0$
since the latter is induced from
the differential on $C$,
and this means that
$$
0 =d^0 d^1 + d^1d^0;
\tag1.6.2.1
$$
it is precisely at this stage where
the sign $(-1)^n$ in (1.3$'$) is needed.
Consequently
$$
d = d^0 + d^1
\tag1.6.2.2
$$
endows 
$\roman{Alt}_A(L,C)$
with the structure of a chain complex over $R$.
\smallskip
Notice that in the special case where $C$ is concentrated in degree zero,
that is, where $C$ is just an  $A$-module,
$Z^0(\roman{End}(C)) =\roman{End}(C)$,
and the differential $d$ boils down to the usual
Lie-algebra cohomology
differential (1.3).
\smallskip
\noindent
{\smc 1.6.3. 
Differential graded
algebras over $A$.}
We shall refer to a differential graded algebra $E$ over $A$ 
as a {\it differential graded\/} $(A,L)$-{\it algebra\/}
if the structure maps
$m \colon E \otimes _A E \to E$
and
$\eta\colon A \to E$
are morphisms of
$(A,L)$-chain complexes.
Equivalently,
write $Z^0(\roman{Der}(E))$
for the homogeneous derivations of $E$ 
of degree zero
that are also
chain maps;
then $E$ is a
differential graded $(A,L)$-algebra
if and only if the structure map
$\omega$
factors through
an action
$\omega \colon L \to Z^0(\roman{Der}(E))$
(denoted still by $\omega$ with a slight abuse of notation).
The same kind of reasoning as above shows that,
for a  differential graded $(A,L)$-algebra $E$,
the pairing (1.5$'$)
and the differential (1.6.2.2)
endow
the $A$-multilinear forms $\roman{Alt}_A(L,E)$
with values in $E$
with the structure of a 
differential
graded
algebra
over $R$.
In fact,
in view of
what was said above, the operator
$d^1$ endows
$\roman{Alt}_A(L,E)$
with the structure of a differential graded algebra,
and it is well known that so does $d^0$
and, 
by virtue of (1.6.2.1),
the operator
(1.6.2.2)
is a differential, too.
\smallskip
Notice in the special case where $E$ is concentrated in degree zero,
that is, where $E$ is just an algebra over $A$,
$Z^0(\roman{Der}(E)) =\roman{Der}(E)$,
and the differential $d$ boils down to the usual
Lie algebra cohomology
differential (1.3).
\smallskip
\noindent
{\smc 1.6.4. 
Differential graded
coalgebras over $A$.}
Recall that a {\it differential graded coalgebra\/} $C$ over $A$ 
is a chain complex $C$ over $A$
together with  chain maps
$\Delta \colon C \to C \otimes_A C$
and
$\varepsilon \colon C \to A$
over $A$
that 
satisfy the usual properties.
Recall also that a {\it coderivation\/}
$\phi \colon C \to C$
of a differential graded coalgebra $C$ is a
morphism of the underlying graded $R$-modules
so that
the diagram
$$
\CD
C
@>{\phi}>>
C
\\
@V{\Delta}VV
@V{\Delta}VV
\\
C\otimes_A C
@>{\phi \otimes_A \roman{Id} + \roman{Id} \otimes_A \phi}>>
C\otimes_A C
\endCD
$$
is commutative.
Notice that the differential of $C$ is itself a coderivation.
We shall refer to a differential graded coalgebra $C$ over $A$ 
 as a {\it differential graded\/} $(A,L)$-{\it coalgebra\/}
if the structure maps
$\Delta \colon C \to C \otimes _A C$
and
$\varepsilon \colon C \to A$
are morphisms of
$(A,L)$-chain complexes.
Equivalently,
write $Z^0(\roman{Coder}(C))$
for the homogeneous coderivations of $C$ 
of degree zero
that are also
chain maps;
then $C$ is a
 differential graded $(A,L)$-coalgebra
if and only if the structure map
$\omega$
factors through
an action
$\omega \colon L \to Z^0(\roman{Coder}(C))$
(denoted still by $\omega$, with a slight abuse of notation).
\smallskip
For later reference
we reproduce the notions of
cofree graded  coalgebra and 
that of
cofree graded cocommutative
coalgebra.
Let $Y$ be a graded $A$-module. 
Consider
the graded tensor coalgebra $(T_A'[Y],\Delta)$ over $A$ on $Y$.
It may be written
$$
T_A'[Y] = \sum { }^{\oplus} T_A^{(n)}[Y],
\tag1.6.4.1
$$
where
$T_A^{(n)}[Y] = Y^{\otimes_A ^n}$,
and where,
for $0 \leq k \leq n$,
the 
component
$$
T_A^{(n)}[Y]
\longrightarrow
T_A^{(k)}[Y]
\otimes_A
T_A^{(n-k)}[Y]
\tag1.6.4.2
$$
of the
diagonal map $\Delta$
is the obvious isomorphism; we note that the direct sum and the tensor product 
are here understood in the graded sense.
For later reference, we spell out the following.

\proclaim{Finiteness hypothesis 1.6.4.3.f}
The graded $A$-module $Y$ is concentrated in non-negative degrees
and
 finitely generated or 
has the property that
$Y_0$ is finitely generated
and $Y_i$ is non-zero only in finitely many degrees $i>0$.
\endproclaim

Let $\pi \colon T_A'[Y] \to Y$
be the canonical projection.
The tensor coalgebra satisfies the following universal property:

\proclaim{1.6.4.3}
Suppose that 
$Y$ 
satisfies the finiteness hypothesis
{\rm (1.6.4.3.f)}.
Then
given any  graded coalgebra $C$ over $A$ and 
a morphism 
$
\phi
\colon
C
\longrightarrow
Y
$
of graded $A$-modules,
there is a unique 
morphism
$
\Phi
\colon
C
\longrightarrow
T_A'[Y]
$
of graded coalgebras over $A$ so that
$\pi \Phi = \phi$.
\endproclaim

Indeed,
given 
$
\phi
\colon
C
\longrightarrow
Y,
$
let
$\Phi = \sum \phi_i$, where
$\phi_0 = \varepsilon$ and where,
for $i \geq 1$,
$\phi_i $
denotes the composite
$$
\phi_i 
\colon
C
@>{\Delta^{(i)}}>>
C^{\otimes_A^i}
@>{\phi}>>
Y^{\otimes_A^i}.
$$
Here $\Delta^{(i)}$
refers to some corresponding iterate of the diagonal map;
it does not matter which one we take,
by coassociativity.
In view of the finiteness hypothesis
(1.6.4.3.f),
the map
$\Phi$ is  well defined,
that is, in each degree, only finitely many terms 
$\phi_i$ are non-zero.
\smallskip
In other words,               
under the finiteness hypothesis (1.6.4.3.f),
the triple $(T_A'[Y],\Delta, \pi)$
constitutes the {\it cofree graded coalgebra\/}
on $Y$ over $A$.
We note that, in case $Y$ does not satisfy a finiteness hypothesis
of the kind (1.6.4.3.f), the tensor coalgebra
on $Y$
will not satisfy the universal
property; it must then be completed, and so must be
the tensor product
which is the target for the diagonal map.
In the application in Section 3 this problem will not occur.
Moreover, 
the 
obvious morphism
$\eta \colon A \to T_A'[Y]$
of coalgebras over $A$
endows the tensor coalgebra
with the structure of a 
{\it coaugmentation\/}.
Consequently
the tensor coalgebra is filtered by the {\it coaugmentation filtration\/}
(see e.~g. \cite\munkholm).
This filtration 
of $T_A'[Y]$
is the same as that 
by the length of tensors
in each  summand $T_A^{(n)}[Y]$.
\smallskip
Recall that
a  graded coalgebra $C$ over $A$
is {\it cocommutative\/} provided
the composite
of the diagonal $\Delta$
with the interchange map
$$
C \otimes _A C \longrightarrow C \otimes_A C,
\quad a \otimes_A b \mapsto (-1)^{|a||b|} b \otimes_A a,
$$
coincides with $\Delta$.  
Recall that
the 
{\it cofree graded cocommutative coalgebra on\/} 
(the graded $A$-module)
$Y$ {\it in the category
of $A$-modules\/} is a pair
$(S_A'\lbrack Y \rbrack, \pi)$, where $S_A'\lbrack Y \rbrack$
is a graded cocommutative  coalgebra over $A$
and where 
$
\pi \colon 
S_A'\lbrack Y \rbrack
\longrightarrow
Y
$
is a morphism of $A$-modules having the
following universal property:

\proclaim{1.6.4.3.c}
Given any  graded cocommutative coalgebra $C$ over $A$ and 
a morphism 
\linebreak
$
\phi
\colon
C
\longrightarrow
Y
$
of $A$-modules,
there is a unique 
morphism
$
\Phi
\colon
C
\longrightarrow
S_A'\lbrack Y \rbrack
$
of graded (commutative) coalgebras over $A$ so that
$\pi \Phi= \phi$.
\endproclaim

Instead of
\lq\lq cofree graded cocommutative coalgebra on $Y$\rq\rq \ 
we shall also say
{\it graded symmetric coalgebra on\/} $Y$
(over $A$).
It is clear that the graded symmetric coalgebra on $Y$
is unique up to isomorphism,
if it exists.
\smallskip
We now reproduce a construction
of the graded symmetric coalgebra on $Y$:
Let $(T_A'[Y],\Delta)$ be the tensor coalgebra over $A$ on $Y$.
For $n \geq 1$,
let $T_A^{(n)}[Y]$ be its homogeneous degree $n$ component,
let
the symmetric group $S_n$ on $n$ letters
act on
$T_A^{(n)}[Y]$
in the graded sense
in the obvious way, 
that is to say, 
for any $y=y_1\otimes_A\dots\otimes_A y_n \in T_A^{(n)}[Y]$ and
for a transposition $\tau = (i,j)$ avec $i<j$, we have
$$
\tau(y_1\otimes_A\dots\otimes_A y_i\otimes_A \dots \otimes_A
 y_j\otimes_A \dots\otimes_A y_n)
=\varepsilon(i,j,y)
(y_1\otimes_A\dots\otimes_A y_j\otimes_A \dots \otimes_A
 y_i\otimes_A \dots\otimes_A y_n)
$$
where
$$
\varepsilon(i,j,y)=
(-1)^{|y_i||y_j| + (|y_{i+1}|+\dots+|y_{j-1}|)(|y_i|+|y_j|)},
$$
and let
$(S_A')^n[Y]$ be the submodule
of $S_n$-invariants.
Let
$S_A'[Y] = \oplus (S_A')^n[Y]$,
and let
$\pi \colon S_A'[Y] \to Y$
be the obvious projection.
For any degree $n$,
given
$0 \leq k \leq n$,
the component (1.6.4.2) of the restriction of the diagonal map
$\Delta$
to
$T_A^{(n)}[Y]$ 
then maps an element
$x \in (S_A')^n[Y]$,
that is, an element
$x \in T_A^{(n)}[Y]$
invariant
under $S_n$,
to an element in
$T_A^{(k)}[Y] \otimes _A T_A^{(n-k)}[Y]$
invariant under
$S_k \times S_{n-k}$,
with respect to the obvious action of
the latter group on
$T_A^{(k)}[Y] \otimes _A T_A^{(n-k)}[Y]$;
indeed, since by construction,
the morphism (1.6.4.2) is
just the obvious isomorphism,
this amounts to the fact that
$x \in T_A^{(n)}[Y]$
is invariant under 
$S_k \times S_{n-k}$,
with respect to the corresponding obvious
embedding
of $S_k \times S_{n-k}$ into $S_n$.
Consequently
the diagonal map $\Delta$
on $T_A'[Y]$ passes to one on
$S_A'[Y]$ which we still denote by $\Delta$. 
Moreover, the 
coaugmentation
map $\eta$
for $T_A'[Y]$
yields a coaugmentation
map $\eta \colon A \to S_A'[Y]$
for
$S_A'[Y]$.

\proclaim{Proposition 1.6.4.4}
Suppose that 
$Y$ 
satisfies the finiteness hypothesis
{\rm (1.6.4.3.f)}.
Then $(S_A'[Y], \Delta,\pi)$
is the graded symmetric coalgebra on $Y$ over $A$.
\endproclaim

\demo{Proof}
It is clear that,
under the finiteness hypothesis (1.6.4.3.f),
$(S_A'[Y], \Delta,\pi)$
satisfies the universal property (1.6.4.3.c). \qed
\enddemo
\smallskip

In view of what was said above,
the coaugmentation filtration 
of $S_A'[Y]$
is the same as that 
by the length of invariant tensors
in $T_A^{(n)}[Y]$.

\proclaim{Proposition 1.6.4.5}
For a graded $(A,L)$-module $Y$,
the
usual formula
$$
\alpha (y_1 \otimes_A y_2 \otimes_A \dots \otimes _A y_n)
=
\sum
y_1 \otimes_A y_2 \otimes_A \dots \otimes _A
\alpha(y_i)\otimes _A \dots  \otimes _A y_n,
\quad \alpha \in L,
$$
endows
\roster
\item
the graded tensor
coalgebra
$(T_A'[Y],\Delta)$
with a structure
of a graded $(A,L)$-coalgebra, and
\item
the graded tensor
algebra
$(T_A[Y],m)$
with a structure
of a graded $(A,L)$-algebra.
\endroster
Furthermore, if $Y$
 satisfies the finiteness hypothesis 
{\rm (1.6.4.3.f)} so that
the graded symmetric coalgebra
$(S_A'[Y],\Delta)$
exists, the latter
inherits
 a structure
of a graded $(A,L)$-coalgebra.
\endproclaim

\demo{Proof} This is left to the reader. \qed \enddemo

When $Y$ is concentrated in odd degrees,
$S'_A[Y]$
is just the {\it (graded) exterior coalgebra\/} $\Lambda'_A[Y]$;
when $Y$ is, furthermore, projective as a graded $A$-module,
as graded $A$-modules,
the exterior algebra $\Lambda_A[Y]$
and coalgebra $\Lambda'_A[Y]$ coincide; the structures
in fact combine to that of a graded Hopf algebra.
More precisely,
the usual diagonal map
$$
\Delta \colon
\Lambda_A[Y]
\longrightarrow
\Lambda_A[Y]\otimes _A\Lambda_A[Y]
\tag1.6.4.6
$$
determined by
$$
\Delta(v) = v \otimes_A 1 + 1 \otimes_A v,\quad v \in Y,
\tag1.6.4.7
$$
endows 
the graded exterior algebra
$\Lambda_A[Y]$ with the structure of
a graded commutative and graded
cocommutative Hopf algebra and hence in particular with that
of a graded cocommutative coalgebra;
see e.~g. {\smc Mac Lane}~\cite\maclaboo\ for details.
The latter is precisely
the graded exterior $A$-coalgebra structure on $Y$.
Since 
the property of being a Hopf algebra implies in particular
that its diagonal map is multiplicative, the rule (1.6.4.7)
in fact completely determines (1.6.4.6).
Explicitly, 
given
$x_1,\dots,x_{p+q} \in Y$,
the value 
$$
\Delta(x_1 \wedge_A\dots\wedge_A x_{p+q}) \in \Lambda_A[Y]
$$ 
is given by the formula
$$
\aligned
\Delta (x_1\wedge_A\dots\wedge_A x_{p+q})&
\\ 
=
\sum_{\sigma}\roman{sign}(\sigma)
&(x_{\sigma(1)}\wedge_A\dots\wedge_A x_{\sigma(p)})
\otimes_A(x_{\sigma(p+1)}\wedge_A\dots\wedge_A x_{\sigma(p+q)}),
\endaligned
\tag1.6.4.8 
$$
where
$\sigma$
runs through $(p,q)$-shuffles
and where $\roman{sign}(\sigma)$ refers to the sign of $\sigma$.
\smallskip
Likewise,
when $Y$ is concentrated in even degrees,
$S'_A[Y]$ is 
the {\it (graded) symmetric coalgebra\/}
$\Sigma'_A[Y]$
in the category of $A$-modules,
but
the relationship
between the graded symmetric algebra 
$\Sigma_A[Y]$
and the graded 
symmetric coalgebra cannot in general be explained in terms of 
an isomorphism of Hopf algebras.
This relationship is actually folk-lore
but difficult to locate in the literature.
We therefore spell out some of the details.
\smallskip
We still suppose that $Y$ is projective as a graded $A$-module.
The 
diagonal map
$$
\Delta
\colon
Y
@>>>
Y \oplus Y,
\quad
\Delta (y) = (y,y),
\ y \in Y,
$$
induces
a diagonal map
$$
\Delta
\colon
\Sigma_A[Y]
@>>>
\Sigma_A[Y]
\otimes_A
\Sigma_A[Y]
\cong
\Sigma_A[Y \oplus Y]
\tag1.6.4.9
$$
which endows $\Sigma_A[Y]$ with the structure of a (graded)
commutative and cocommutative
Hopf algebra.
Explicitly, 
given
$x_1,\dots,x_{p+q} \in Y$,
the value 
$$
\Delta(x_1 x_2 \dots x_{p+q}) \in \Sigma_A[Y]
$$ 
is given by the formula
$$
\Delta (x_1x_2\dots x_{p+q}) 
=
\sum_{\sigma}
(x_{\sigma(1)}\dots x_{\sigma(p)})
\otimes_A(x_{\sigma(p+1)}\dots x_{\sigma(p+q)})
\tag1.6.4.10 
$$
where
$\sigma$
runs through $(p,q)$-shuffles.
This diagonal map is referred to as
{\it shuffle coproduct\/}.
In view of the universal property (1.6.4.3.c)
of the symmetric
coalgebra
$\Sigma'_A[Y]$,
the canonical projection
$\phi \colon 
\Sigma_A[Y] \to Y$
induces a morphism
$$
\Sigma_A[Y]
@>>>
\Sigma'_A[Y]
\tag1.6.4.11
$$
of graded (commutative)
coalgebras.
Furthermore,
again in view of the universal property (1.6.4.3.c)
of the symmetric
coalgebra
$\Sigma'_A[Y]$,
addition
$$
Y \oplus Y
@>>>
Y,
\quad
(y_1,y_2)
\mapsto y_1+y_2,
\quad y_1, y_2 \in Y,
$$
induces a multiplication map
$$
\Sigma'_A[Y]
\otimes_A
\Sigma'_A[Y]
\cong
\Sigma'_A[Y \oplus Y]
@>>>
\Sigma'_A[Y]
\tag1.6.4.12
$$
which endows $\Sigma'_A[Y]$ with the structure of a (graded)
commutative and cocommutative
Hopf algebra as well,
and
(1.6.4.11)
is a morphism of Hopf algebras.
Since
$\Sigma_A[Y]$
is the free graded commutative algebra on $Y$,
the  morphism of
graded algebras which underlies
(1.6.4.11)
may also be obtained as the unique morphism
of graded algebras
induced by the canonical
inclusion
from $Y$ into 
$\Sigma'_A[Y]$,
viewed as a graded algebra.
\smallskip
When we dualize 
$\Sigma'_A[Y]$
and
$\Sigma_A[Y]$,
we obtain the graded Hopf algebras
$$
\Sigma_A[Y^*]
=
\roman{Hom}_A(\Sigma'_A[Y],A)
$$
and
$$
\Sigma^*_A[Y]
=
\roman{Hom}_A(\Sigma_A[Y],A)
=
\Sigma'_A[Y^*]
$$
respectively,
where
$Y^*= \roman{Hom}_A(Y,A)$;
here
$\Sigma_A[Y^*]$
is the algebra of {\it graded polynomial functions\/}
on $Y$
and
$\Sigma^*_A[Y]$
that of {\it graded symmetric functions\/}
on $Y$;
the multiplication
on the algebra of graded symmetric functions
is the dual of the shuffle coproduct
(1.6.4.9)
and hence the usual {\it shuffle product\/}.
The dual
of
(1.6.4.11)
is the canonical map
$$
\Sigma_A[Y^*]
@>>>
\Sigma^*_A[Y]
\tag1.6.4.13
$$
from the Hopf algebra of  graded polynomial functions
to that of graded symmetric functions
on $Y$.
It 
sends a polynomial function to the corresponding
symmetric function and
is formally exactly of the same kind as
(1.6.4.11),
with $Y^*$ instead of $Y$.
\smallskip
When $Y$ is $A$-free, after a choice of basis
$\{e_1,e_2,\dots\}$ of $Y$ has been made
and when
$\{\xi_1,\xi_2,\dots\}$ refers to the dual basis of $Y^*$,
$\Sigma_A[Y^*]$  
is the polynomial $A$-Hopf algebra 
$A[\xi_1,\xi_2,\dots]$
on
$\xi_1,\xi_2,\dots $
whereas
$\Sigma^*_A[Y]$
is the
{\it divided polynomial $A$-Hopf algebra\/}
$\Gamma [\xi_1,\xi_2,\dots ]$
on
$\xi_1,\xi_2,\dots $,

that is, 
$\Sigma^*_A[Y]$ is
the graded commutative $A$-Hopf algebra
generated by
$\gamma_k\xi_j$,
$k,j \geq 1$,
subject to the relations
$$
k! 
\gamma_k\xi_j = \xi_j^k,
\quad
k,j \geq 1,
$$
with $A$-coalgebra structure $\Delta$ determined by
$$
\Delta
\gamma_k\xi_j
=
\sum_{u+v = k}
\gamma_u\xi_j
\otimes_A
\gamma_v\xi_j.
$$
See \cite\cartanse\ for more details on divided powers.
The map
(1.6.4.13)
is then the obvious one
which sends 
the multiplicative generator $\xi_j$
to the multiplicative generator $\gamma_1\xi_j = \xi_j$
but 
(1.6.4.13) is not
in general an isomorphism.
In characteristic zero
it is an isomorphism, though, since
we can then define the divided power operations
in $A[\xi_1,\xi_2,\dots]$
by
$$
\gamma_k\xi_j = \frac 1{k!}\xi_j^k,
\quad
j,k \geq 1;
$$
the inverse mapping of
(1.6.4.13) is then usually referred to
as {\it polarization\/}.
\smallskip
\noindent
{\smc 1.6.5. 
Cup products.}
Let
$C$ 
be a differential graded coalgebra
in the category of $A$-modules,
with structure map
$C
@>{\Delta}>>
C\otimes_A C$,
and let
$$
\mu_A
\colon
M'
\otimes_A
M''
\longrightarrow
M
\tag1.6.5.1
$$
be a pairing
of differential graded $A$-modules, that is, of
chain complexes in the category of $A$-modules.
Given morphisms $a\colon C\to M'$ and $b\colon C\to M''$ 
(say) of the underlying 
graded modules, the
{\it cup product} of $a$ and $b$
with respect to $\mu_A$, written $a\cup b$, is the composite
$$
C
@>{\Delta}>>
C\otimes_A C
@>{a\otimes_A b}>>
M'
\otimes_A
M''
@>{\mu_A}>>
M;
\tag1.6.5.2
$$
see
\cite{\perturba,\,\husmosta,\,\maybook,\,\munkholm}. 
Here the tensor product
$a\otimes_A b$ is the {\it graded one\/}, that is to say,
$$
a\otimes_A b(x\otimes_A y) = (-1)^{|b||x|} \mu_A(a(x) \otimes_A b(y)),
\quad x \in M',\ y \in M''.
\tag1.6.5.3
$$
The resulting pairing
$$
\cup \colon
\roman{Hom}_A(C,M')
\otimes_A
\roman{Hom}_A(C,M'')
@>>>
\roman{Hom}_A(C,M),
\quad
(a,b)\mapsto a\cup b,
\tag1.6.5.4
$$
is the usual
{\it cup pairing\/}
with respect to $\mu_A$; it 
is associative in the obvious way.
\smallskip
To relate this kind of cup pairing with
the pairing (1.5$''$),
given an $(R,A)$-Lie algebra $L$, let
$Y = sL$, the {\it suspension\/}
$sL$ of $L$; for the present purposes this means that
$sL$ is just $L$ except that its elements are regraded up by one.
Given a pairing 
 of differential graded $A$-modules 
of the kind (1.6.5.1),
we then have
the corresponding pairing
(1.6.5.4).
In the special case where
the pairing (1.6.5.1)
is one of ungraded $A$-modules, viewed
as differential graded $A$-modules
concentrated in degree zero,
the resulting pairing
(1.6.5.4)
is exactly the same as the pairing (1.5$''$) above.
It is exactly at this stage where
the sign
$(-1)^{|\alpha||\beta|}$
in the formula (1.4.2) above is needed.
\smallskip
In the general (graded) case, for an arbitrary
differential graded coalgebra $C$
and  a differential graded  algebra $U$,
both
in the category of $A$-modules,
the cup
product turns $\roman{Hom}_A(C,U)$ into a differential graded algebra 
in the category of $A$-modules,
with
unit $\eta\varepsilon$ and
differential $D$ given by
$$
Df  = df + (-1)^{(|f|+1)}fd,
\tag1.6.5.5
$$
for homogeneous $f \colon C \to U$;
further, if
$C$ and $U$ have a coaugmentation map $\eta$ and
augmentation map $\varepsilon$, respectively, the assignment
$\varphi\mapsto\varepsilon\varphi\eta$ 
yields
an augmentation map for this differential graded algebra.

\beginsection 2. Extensions

The algebraic analog
of 
an \lq\lq Atiyah sequence\rq\rq\ 
or of a
\lq\lq transitive Lie algebroid\rq\rq,
see e.~g.~\cite\mackbook\ or (2.2) below,
is an extension of Lie-Rinehart algebras.
In the present section
we study
such extensions
by means of generalizations of the
usual notions of
connection
and curvature
in a principal bundle.
\smallskip
Let $L'$, $L$, $L''$ be $(R,A)$-Lie algebras.
An
{\it extension\/} of
$(R,A)$-Lie algebras is a short exact sequence
$$
\roman {\bold e}
\colon
0
@>>>
L'
@>>>
L
@>p>>
L''
@>>>
0
\tag2.1
$$
in the category of $(R,A)$-Lie algebras;
notice in particular that the
Lie algebra $L'$ necessarily acts trivially on $A$.
If also
$
\bar
{\roman {\bold e}}
\colon
0
\to
L'
\to
\bar L
\to
L''
\to
0
$
is an extension of
$(R,A)$-Lie algebras,
as usual, $\roman {\bold e}$ and $\bar {\roman {\bold e}}$ are
said to be {\it congruent\/},
if there is a
morphism
$
(\roman {Id},\cdot,\roman {Id})
\colon \roman {\bold e} \longrightarrow\bar {\roman {\bold e}}
$
of extensions of
$(R,A)$-Lie algebras.

\smallskip
\noindent
{\smc Remark 2.2.}
Let $N$ be a smooth finite dimensional manifold,
let $A$ be the algebra of smooth functions on 
$N$, 
and let
$\xi \colon P \to N$ be a principal bundle,
with structure group $G$ 
acting from the right.
The
vertical subbundle
$\psi\colon V \to P$ 
of the tangent bundle $\tau_P$
of $P$
is well known to be
trivial,
having as fibre
the Lie algebra $\fra g$ of $G$,
that is,
$V \cong P \times \fra g$.
Dividing out the actions of $G$ from the right, we
obtain
an extension
$$
0
@>>>
\roman{ad}(\xi)
@>>>
\tau_P/G
@>>>
\tau_N
@>>>
0
\tag2.2.1
$$
of vector bundles over $N$, 
where
$\tau_N$ is the tangent bundle of $N$.
This sequence has been introduced
by {\smc Atiyah}~\cite\atiyaone\ (Theorem 1)
and is now usually called
the {\it Atiyah sequence \/} of the principal 
bundle $\xi$;
here $\roman{ad}(\xi)$
is the bundle associated to the principal bundle
by the adjoint representation of $G$ on its Lie algebra $\fra g$.
A complete account to Atiyah sequences
may be found in
App. A of
\cite\mackbook.
The spaces
${\fra g}(\xi)= \Gamma\roman{ad}(\xi)$ and
$E(\xi)=\Gamma(\tau_P/G)$
of sections
inherit obvious structures of
Lie algebras, 
in fact of $(\Bobb R,A)$-Lie algebras, and
$$
0
@>>>
\fra g(\xi)
@>>>
E(\xi)
@>>>
\roman{Vect}(N)
@>>>
0
\tag2.2.2
$$
is an extension of  $(\Bobb R,A)$-Lie algebras;
here $\roman{Vect}(N)=\Gamma(\tau_N)$,
the Lie algebra of vector fields on $N$, and
$\fra g(\xi)$
is in an obvious way the Lie algebra
of the {\it group\/}
of {\it gauge transformations\/}
of $\xi$.
\smallskip
We now generalize the classical notions of
principal connection and curvature:
Let $\roman {\bold e}$ be an extension (2.1) of $(R,A)$-Lie algebras,
and suppose that it splits in the category of $A$-modules;
this will e.~g. hold if
$L''$ is projective as an $A$-module.
Then $\roman {\bold e}$ may be represented
by a 2-cocycle: Let
$
\omega 
\colon
L''
\to
L
$
be a section of $A$-modules
for the projection $p \colon L \to L''$.
We refer
to $\omega$ as an $\roman {\bold e}$-{\it connection\/}.
Given an
$\roman {\bold e}$-connection, define
the corresponding
($\roman {\bold e}$-){\it curvature\/}
$
\Omega
\colon
L'' \otimes_A L''
\to
L'
$
as the morphism $\Omega$ of $A$-modules
satisfying
$$
[\omega(\alpha),\omega(\beta)]
=
\omega[\alpha,\beta]
+ \Omega(\alpha,\beta)
\tag2.3
$$
for every $\alpha, \beta  \in L''$.
The usual reasoning reveals that $\Omega$ is indeed well defined
as an alternating $A$-bilinear 2-form on $L''$ with values in $L'$;
under the circumstances
of (2.2),
this amounts to
$\Omega$ being a tensor.
These notions of
$\roman {\bold e}$-connection
and
 $\roman {\bold e}$-curvature
generalize
the concepts of  principal connection
and  principal curvature;
indeed,
under the circumstances
of (2.2),
they come down to their descriptions
in the language of Atiyah sequences
due to {\smc Mackenzie}~\cite\mackbook.
\smallskip
In~\cite\poiscoho\ we have shown that
in view of  the Jacobi identity  in $L$
the morphism $\Omega$ must satisfy a
2-cocycle condition
phrased in terms of a suitable notion of
covariant derivative
which generalizes the usual Bianchi identity.
We now explain this somewhat more formally:
\smallskip
The adjoint representation
$\roman{ad} \colon L \to \roman{End}(L')$
endows
the $A$-Lie algebra  $L'$ 
with a structure of an
$(A,L)$-module, in fact with that of an
$(A,L)$-Lie algebra
$$
\roman{ad} \colon L \longrightarrow \roman{Der}(L')
\tag2.4
$$
in the sense
of (1.6.1).
In particular, we have the
chain complex
$\roman{Alt}_A(L,L')$
with the 
differential $d_L$ given by (1.3).
Furthermore, the projection $p \colon L \to L''$ induces an injection
$$
p^*
\colon
\roman{Alt}_A(L'',L')
\longrightarrow
\roman{Alt}_A(L,L')
$$
of graded $A$-modules,
and the chosen $\roman {\bold e}$-connection $\omega$
induces a
surjection
$$
\omega^*
\colon
\roman{Alt}_A(L,L')
\longrightarrow
\roman{Alt}_A(L'',L')
$$
of graded $A$-modules.
Define the operator
$D^{\omega}$
of {\it covariant derivative\/}
on the graded $A$-module
$\roman{Alt}_A(L'',L')$
as the composite
$$
D^{\omega} = \omega^* d_Lp^*  \colon
\roman{Alt}_A(L'',L')
\longrightarrow
\roman{Alt}_A(L'',L').
\tag2.5.1
$$  
When we write  out this operator, we obtain the usual formula
$$
\alignedat 1
(D^{\omega}f)&(\alpha_1,\dots,\alpha_n)
\\
&=
\quad 
(-1)^n
\sum_{i=1}^n (-1)^{(i-1)}
\roman {ad}(\omega(\alpha_i))
(f (\alpha_1, \dots\widehat{\alpha_i}\dots, \alpha_n))
\\
&\phantom{=}+\quad
(-1)^n
\sum_{j<k} (-1)^{(j+k)}f(\lbrack \alpha_j,\alpha_k \rbrack,
\alpha_1, \dots\widehat{\alpha_j}\dots\widehat{\alpha_k}\dots,\alpha_n).
\endalignedat
\tag2.5.2
$$
It is readily seen that the Jacobi identity in $L$ boils down to
the identity
$$
D^{\omega}(\Omega) = 0
\tag2.5.3
$$
which 
generalizes the Bianchi identity;
we refer to (2.5.3) as  the {\it generalized Bianchi identity\/}.
\smallskip
The 2-cocycle $\Omega$ is uniquely determined
by 
$\roman {\bold e}$
up to a coboundary;
see Section 2 in~\cite\poiscoho\ for details.
Moreover, from that paper, we recall the following.

\proclaim{Theorem 2.6}
Given an $(R,A)$-Lie algebra
$L''$
and  an $(A,L'')$-module
$L'$,
viewed as an abelian
$(A,L'')$-Lie algebra,
the assignment, 
to an extension
which splits in the category of $A$-modules and
realizes the $(A,L'')$-module structure on $L'$,
of its 2-cocycle $\Omega \in \roman{Alt}_A(L'',L')$,
yields  a bijective correspondence between the congruence classes 
of such extensions of $L'$ by $L''$ 
and the elements of
$\roman H^2(\roman{Alt}_A(L'',L'))$.
\endproclaim

\smallskip
When $L'$ is non-abelian,
the 
generalized Bianchi identity (2.5.3) says that
the 2-form $\Omega$ is a non-abelian 2-cocycle, and hence
$\Omega$
does {\it not\/} lead to a 
cohomology class in 
a naive way.
In this case, the classical argument due to 
{\smc Eilenberg-Mac Lane}~\cite\eilmacon,
see e.~g. (IV.8.8) in {\smc Mac Lane}~\cite\maclaboo,
suitably rephrased for the present case,
shows that
the cohomology group
$\roman H^2_A(\roman{Alt}_A(L'',Z))$
acts faithfully and transitively
on the congruence classes
of extensions of 
$L''$ by $L'$
with the same \lq\lq outer action\rq\rq\ 
of $L''$ on $L'$
where
$Z$ refers to the center of $L'$ as an $A$-Lie algebra.
It seems worthwhile giving some of the details:
At first,
the term \lq\lq outer action\rq\rq\ 
means the following:
The Lie algebra $L'$ acts on itself by means of the
adjoint representation $\roman{ad}\colon L' \to \roman{Der}(L')$
and the image is well  known to be a Lie ideal;
hence the morphism  $\roman{ad}$  admits a cokernel, the Lie algebra
$\roman{ODer}(L')$ of {\it outer\/} derivations of $L'$.
We refer to an arbitrary morphism 
$L'' \to \roman{ODer}(L')$ 
of $R$-Lie algebras as
an {\it outer action\/} of $L''$ on $L'$.
An outer action induces an action
$$
L'' \to \roman{Der}(Z)
$$
that endows the center $Z$ of $L'$ with a structure of
an $(A,L'')$-module, in fact, with that of 
an $(A,L'')$-Lie algebra in the sense of (1.6.1),
with trivial Lie structure on $Z$ understood.
In particular,
the chain complex $\roman{Alt}_A(L'',Z)$
is well defined.
Furthermore,
given an extension $\roman {\bold e}$ of Lie-Rinehart algebras
of the kind {\rm (2.1)},
the corresponding action  {\rm (2.4)}
induces an outer action 
$$
L'' \to \roman{ODer}(L')
\tag2.7.1
$$
of $L''$ on $L'$.
Next, addition induces an operation
$$
+ \colon L' \oplus Z \longrightarrow L'
$$
of $A$-modules which, in turn, induces an operation
$$
+ \colon\roman{Alt}_A(L'',L')
\oplus \roman{Alt}_A(L'',Z )
\longrightarrow
\roman{Alt}_A(L'',Z).
\tag2.7.2
$$
A choice of $\roman {\bold e}$-connection
corresponds to a decomposition of
$L$ as a direct sum $L' \oplus L''$ as $A$-modules, and
the aforementioned argument due to 
{\smc Eilenberg-Mac Lane}~\cite\eilmacon,
suitably rephrased, establishes a proof of the following.

\proclaim{Theorem 2.7}
Given an extension
$\roman {\bold e}$ of Lie-Rinehart algebras
of the kind {\rm (2.1)},
the operation 
{\rm (2.7.2)}
induces
a faithful and transitive action of
the cohomology group
$\roman H^2(\roman{Alt}_A(L'',Z))$
on the congruence classes
of extensions of 
$L''$ by $L'$
with \lq\lq outer action\rq\rq\ 
{\rm (2.7.1)}
of $L''$ on $L'$
in such a way that, when $\Omega$ is the curvature corresponding to
an extension $\roman {\bold e}$
and when $\rho \in \roman{Alt}^2_A(L'', Z)$
is a 2-cocycle,
the 2-form
$$
\Omega + \rho
$$ 
is the curvature corresponding to
an extension $\roman {\bold e}_{\rho}$
representing the congruence class
obtained by operation upon 
the class of $\roman {\bold e}$
with $[\rho] \in \roman H^2(\roman{Alt}_A(L'',Z))$. \qed
\endproclaim

A version of this result in the framework of
transitive Lie algebroids
may be found in (IV.3.31) of \cite\mackbook.
An $\roman {\bold e}$-connection is said to be {\it flat\/}
if its curvature is zero.
It is clear that an extension $\roman {\bold e}$ admits
a flat
$\roman {\bold e}$-connection
if and only if 
it splits in the category of
$(R,A)$-Lie algebras.
\smallskip
In the next section we shall need a suitable notion of an operator
of covariant derivative
for $(A,L)$-modules. We now explain this briefly:
Let $M$ be an $(A,L)$-module, possibly graded.
Then the restriction 
to $L'$
of the structure
map from $L$ to $\roman{End}_R(M)$
is an action $\phi$ of $L'$ on $M$ in the usual sense
of Lie algebra actions in the category of $A$-modules.
Consider
the chain complex
$\roman{Alt}_A(L,M)$,
with the Lie algebra cohomology
differential $d_L$, cf. (1.3).
The projection $p \colon L \to L''$ induces an injection
$$
p^*
\colon
\roman{Alt}_A(L'',M)
\longrightarrow
\roman{Alt}_A(L,M)
$$
of graded $A$-modules.
Pick an $\roman {\bold e}$-connection $\omega$;
this induces a
surjection
$$
\omega^*
\colon
\roman{Alt}_A(L,M)
\longrightarrow
\roman{Alt}_A(L'',M)
$$
of graded $A$-modules
and hence
an operator
$D^{\omega}$
of {\it covariant derivative\/}
on the graded $A$-module
$\roman{Alt}_A(L'',M)$, given
as the composite
$$
D^{\omega} = \omega^* d_Lp^*  \colon
\roman{Alt}_A(L'',M)
\longrightarrow
\roman{Alt}_A(L'',M).
\tag2.8.1
$$
A version of this may be found in (IV.3.9)
of \cite\mackbook.
When we write out (2.8.1) we obtain the usual formula
$$
\alignedat 1
(D^{\omega}f)&(\alpha_1,\dots,\alpha_n)
\\
&=
\quad 
(-1)^n
\sum_{i=1}^n (-1)^{(i-1)}
(\omega(\alpha_i))
(f (\alpha_1, \dots\widehat{\alpha_i}\dots, \alpha_n))
\\
&\phantom{=}+\quad
(-1)^n
\sum_{j<k} (-1)^{(j+k)}f(\lbrack \alpha_j,\alpha_k \rbrack,
\alpha_1, \dots\widehat{\alpha_j}\dots\widehat{\alpha_k}\dots,\alpha_n).
\endalignedat
\tag2.8.2
$$
It is manifest that
the 
$\roman {\bold e}$-curvature
$
\Omega
\colon
L'' \otimes_A L'' \to L'
$
of $\omega$,
combined with the
Lie algebra action
$\phi \colon
L'
\to
\roman{End}_A(M)
$
of $L'$
on $M$ in the category of $A$-modules
mentioned earlier
is then the adjoint of
$$
D^{\omega} D^{\omega}
\colon
M
\longrightarrow
\roman{Alt}^2_A(L'',M).
\tag2.9
$$
Explicitly, with $D = D^{\omega}$,
this is formally the usual formula
$$
\phi\Omega(\alpha,\beta) = 
D_\alpha D_\beta - D_\beta D_\alpha - D_{[\alpha,\beta]}
$$
where
$$
D_{\alpha}(m) = (\omega(\alpha))(m),
\quad \alpha \in L,\ m \in M.
$$
\smallskip
For a general $(R,A)$-Lie algebra $L$
and an $A$-module $M$,
there are notions
of $L$-connection and $L$-curvature,
and these can be realized by an
action on $M$ of a
suitable extension
$E$ of $L$ by $\roman{End}_A(M)$
in the category
of $(R,A)$-Lie algebras.
We do not need the details here;
see e.~g. \cite\poiscoho\ (Section 2).

\beginsection 3. The Chern-Weil construction

As before we suppose that
the extension (2.1) splits in the category of $A$-modules.
Let
$s^2L'$
be the  double suspension
of $L'$, i.~e.
$s^2L'$
is just $L'$,
except that its elements are regraded up by 2.
The graded algebra 
$\roman{Hom}_A(\Sigma'_A[s^2L'],A)^L$
being equipped with the zero differential,
we shall construct a morphism
$$
\roman{Hom}_A(\Sigma'_A[s^2L'],A)^L
\longrightarrow
\roman{Alt}_A(L'',A)
$$
of
differential graded commutative $R$-algebras
whose induced morphism
on homology depends only on the 
congruence class 
(cf. Section 2)
of the
extension (2.1).
When $L'$ is finitely generated and projective,
the graded $A$-algebra
$\roman{Hom}_A(\Sigma'_A[s^2(L')],A)$
may be identified with
the symmetric $A$-algebra 
$\Sigma_A[s^2(L')^*]$
on the $A$-dual
$s^2(L')^*$
of $s^2(L')$ as indicated.
In particular, when $L'$ is $A$-free of finite type,
$\roman{Hom}_A(\Sigma'_A[s^2(L')],A)$
is the polynomial algebra over $A$ on an $A$-basis
of
$\roman{Hom}_A(s^2(L'),A)$.
\smallskip
By assumption,
$L'$ is 
a Lie algebra over $A$ in the usual sense.
Furthermore,
the extension (2.1) splits in the category of
$A$-modules;
let $\omega \colon L'' \to L$ be an $\roman {\bold e}$-connection for
(2.1), and let
$$
\Omega \colon L'' \otimes _A L'' \longrightarrow L'
\tag3.1
$$
be its curvature.
Since $\Omega$ is an $A$-bilinear alternating form, it passes through
the second exterior power
$\Lambda_A^2[L'']$
and hence induces 
a morphism
$$
\Lambda_A^2\Omega
\colon
\Lambda_A^2[L'']
\longrightarrow
L'
$$
of $A$-modules.
Thus, in view of (1.6.4.3.c),
$\Lambda'_A [sL'']$ being endowed with the
graded exterior $A$-coalgebra structure (1.6.4.6),
the curvature $\Omega$ induces a morphism
$$
\Omega_{\sharp}
\colon
\Lambda'_A [sL'']
\longrightarrow
\Sigma'_A[s^2L']
\tag3.2
$$
of graded coalgebras over $A$
in the following way:
Let $\pi \colon \Sigma'_A[s^2L'] \to s^2L'$
be the projection map which is part of the structure of
the graded symmetric coalgebra over $A$
(cf. (1.6.4)),
and let
$$
\Omega_{\natural}
\colon
\Lambda_A^2[sL'']
@>>>
s^2L',
\qquad
\Omega_{\flat}
\colon
\Lambda'_A [sL'']
\longrightarrow
s^2L'
$$
be the homogeneous
degree zero morphisms of graded $A$-modules determined by the requirement
that the diagram
$$
\CD
\Lambda'_A [sL'']
@>{\Omega_{\flat}}>>
s^2L'
\\
@V{\roman{proj}}VV
@V{\roman{Id}}VV
\\
\Lambda_A^2[sL'']
@>{\Omega_{\natural}}>>
s^2L'
\\
@A{\Lambda_A^2s}AA
@A{s^2}AA
\\
\Lambda_A^2[L'']
@>{\Lambda_A^2\Omega}>>
L'
\endCD
\tag3.3
$$
be commutative;
here \lq\lq proj\rq\rq\ refers to the obvious projection
map
from $\Lambda'_A [sL'']$ to $\Lambda_A^2[sL'']$.
We note that
$\Omega_{\natural}$
is just the appropriate rewrite
of
$\Lambda_A^2\Omega$
in the formally appropriate graded setting.
In view of the universal property (1.6.4.3.c)
of the graded symmetric coalgebra,
the induced morphism (3.2)
is determined by
the requirement that
$$
\pi\Omega_{\sharp}= \Omega_{\flat}.
\tag3.4
$$
We refer to
$\Omega_{\sharp}$
as a {\it classifying map for the extension\/} (2.1).
It may be viewed as an algebraic analogue
of the more usual notion of classifying map
in topology and differential geometry.
Since $\Omega$ is unique up to a coboundary,
this classifying map is uniquely determined by
(2.1) up to a non-abelian coboundary in a suitable sense.
The construction
of 
$\Omega_{\sharp}$
is completely formal
and does {\it not\/}
require that $L'$ satisfy any finiteness assumption;
the reason is that we work with the  graded
symmetric coalgebra
$\Sigma'_A[s^2L']$
which,
apart from other formal advantages,
in particular 
removes
the existence problem 
for the cofree cocommutative coalgebra
we would be faced with in general if we had tried an ungraded
construction.
\smallskip
The classifying map (3.2) induces a morphism
$$
\roman{Hom}_A(\Sigma'_A[s^2L'],A)
\longrightarrow
\roman{Alt}_A(L'',A)
\tag3.5
$$
of graded $A$-algebras.
To manufacture a chain map from it,
we observe that 
the adjoint representation
of $L$ on itself induces an action of $L$ on
$L'$
that endows $L'$  with the structure
of an $(A,L)$-module
(in fact, with that of an $(A,L)$-Lie algebra, cf. (1.6.1)).
In view of (1.6.4.5),
this $L$-action on $L'$
induces an
action
$$
\omega_{\roman {\bold e}} \colon L \to \roman{Coder}^0_R(\Sigma'_A[s^2L'])
\tag3.6
$$
which endows the latter with the structure of a graded $(A,L)$-coalgebra
in a sense explained in (1.6.4).
These structures, in turn, induce on
$\roman{Hom}_A(\Sigma'_A[s^2L'],A)$
the structure of a graded commutative $(A,L)$-algebra
in 
a sense explained in (1.6.3).
As usual, for 
$\zeta \colon 
\Sigma'_A[s^2L'] 
\to
A$ and $\alpha \in L$, the result of acting upon
$\zeta$ with $\alpha$ is given by
$$
\alpha(\zeta) =
\alpha \circ \zeta - \zeta\circ \alpha
\colon 
\Sigma'_A[s^2L'] 
\longrightarrow
A
\tag3.7
$$
where, with an abuse of notation,
the operators on $\Sigma'_A[s^2L']$ and $A$ induced by $\alpha$
are denoted by $\alpha$ as well.

\proclaim{Theorem 3.8}
Given an extension 
$\roman {\bold e}$
of $(R,A)$-Lie algebras
of the kind
{\rm (2.1)}
and an $\roman {\bold e}$-connection $\omega$
with curvature
{\rm (3.1)},
the restriction of
{\rm (3.5)}
to the invariants
$\roman{Hom}_A(\Sigma'_A[s^2L'],A)^L$
goes into the cycles
of
$\roman{Alt}_A(L'',A)$.
In other words,
when
$\roman{Hom}_A(\Sigma'_A[s^2L'],A)^L$
is endowed with the zero differential,
the restriction of {\rm (3.5)}
yields a morphism 
$$
\roman{Hom}_A(\Sigma'_A[s^2L'],A)^L
\longrightarrow
\roman{Alt}_A(L'',A)
\tag3.8.1
$$
of
differential graded commutative $R$-algebras.
Furthermore,
the induced morphism
$$
\roman{Hom}_A(\Sigma'_A[s^2L'],A)^L
\longrightarrow
\roman H^{2*}(\roman{Alt}_A(L'',A))
\tag3.8.2
$$
depends only on the 
congruence class of the
extension
{\rm (2.1)}
and not on a particular choice of 
an $\roman {\bold e}$-connection $\omega$ etc.
\endproclaim

We shall refer to the morphism (3.8.2)
as the {\it Chern-Weil 
map\/}
for the extension (2.1)
of $(R,A)$-Lie algebras.
The $R$-algebra
$\roman{Hom}(\Sigma_A'[s^2L'],A)^L$
will be referred to as the 
{\it algebra of characteristic classes for\/}
the extension (2.1).
\smallskip
The proof of (3.8) will be subdivided
into several steps.
It will be convenient to view the morphism
$\Omega_{\sharp}$
of graded coalgebras over $A$ as an element of
the  graded $A$-module
$$
\roman{Hom}_A(\Lambda'_A [sL''],\Sigma'_A[s^2L'])
=
\roman{Alt}_A(L'',\Sigma'_A[s^2L']).
\tag3.9
$$
As already pointed out,
the adjoint action (3.6)
induces an action
of $L$ on
$\Sigma'_A[s^2L']$
that endows the latter with the structure of 
an $(A,L)$-module,
in fact, with that of an
$(A,L)$-coalgebra,
and what is said at the end
of Sections 2 applies,
with $M =\Sigma'_A[s^2L']$.
Write
$$
D^{\omega} \colon
\roman{Alt}_A(L'',\Sigma'_A[s^2L'])
\longrightarrow
\roman{Alt}_A(L'',\Sigma'_A[s^2L'])
\tag3.10
$$
for the corresponding operator of covariant derivative,
cf. (2.8.1).

\proclaim{Lemma 3.11}
The
classifying map
$\Omega_{\sharp}$
satisfies
$$
D^{\omega}(\Omega_{\sharp}) = 0.
\tag3.11.1
$$
\endproclaim

The proof requires some preparations; the proof itself will be given
after (3.16) below.
At first we
recall from (1.6.4) that,
by construction,
$\Sigma'_A[s^2L'] \subseteq T'_A[s^2L']$;
we denote the inclusion by $\iota$.
We recall that the graded module underlying the graded tensor coalgebra
$T_A'[s^2L']$ over $A$
coincides with that underlying
the graded tensor algebra
$T_A[s^2L']$ over $A$
(but beware, the algebra and coalgebra structures are {\it not\/} compatible).
The morphism $\iota$ induces an injection
$$
\roman{Hom}_A(\Lambda'_A [sL''],\Sigma'_A[s^2L'])
\longrightarrow
\roman{Hom}_A(\Lambda'_A [sL''],T_A[s^2L'])
\tag3.12
$$
of graded $A$-modules.
Moreover,
in view of
(1.6.5),
the adjoint action 
of $L$ on $L'$
induces also an  action
$$
\omega_{\roman {\bold e}} \colon L \longrightarrow \roman{Der}(T_A[s^2L'])
\tag3.13
$$
of $L$ on $T_A[s^2L']$
which endows
$T_A[s^2L']$
with the structure
of an
$(A,L)$-algebra,
and we can
apply what is said at the end of Section 2,
with $M=T_A[s^2L']$:
Write
$$
D^{\omega} \colon
\roman{Alt}_A(L'',T_A[s^2L'])
\longrightarrow
\roman{Alt}_A(L'',T_A[s^2L'])
\tag3.14
$$
for the corresponding operator of covariant derivative.
Since the values of 
$\omega_{\roman {\bold e}}$
and hence those of
the composite
$
\omega_{\roman {\bold e}}\omega
\colon 
L''
@>>>
\roman{Der}(T_A[s^2L'])
$
lie in the graded $R$-module
$\roman{Der}(T_A[s^2L'])$
of derivations rather than
in the  full
graded $R$-module
$\roman{End}(T_A[s^2L'])$
of $R$-linear endomorphisms,
$D^{\omega}$
is a {\it derivation\/} 
over the ground ring $R$
of the graded $A$-algebra
$$
\roman{Hom}_A(\Lambda'_A [sL''],T_A[s^2L'])
=\roman{Alt}_A(L'',T_A[s^2L']),
$$
the $A$-algebra structure being the shuffle product,
or cup product, cf. (1.6.4) and
(1.6.5), with respect to the
graded exterior $A$-coalgebra structure
on $\Lambda'_A [sL'']$
and
graded $A$-algebra structure on $T_A[s^2L']$.
Finally, the injection (3.12)
is manifestly compatible with the
operations (3.10) and (3.14)
of covariant derivative.

\proclaim{Lemma 3.15}
With respect to the graded 
$A$-coalgebra and
$A$-algebra structures on
$\Lambda_A' [sL'']$ and $T_A[s^2L']$,
respectively,
and with the notation $\Omega_{\flat}$ introduced in
{\rm (3.3)},
the morphism
$$
\iota \Omega_{\sharp}
\colon
\Lambda'_A [sL'']
\longrightarrow
T'_A[s^2L']
$$
may be written
$$
\iota \Omega_{\sharp}
=
\roman{Id}_A +
\Omega_{\flat}  
+
(\Omega_{\flat})\cup (\Omega_{\flat})
+ (\Omega_{\flat})\cup(\Omega_{\flat})\cup (\Omega_{\flat})
+ \dots,
\tag3.15.1
$$
where we do not distinguish in notation between
$$
\Omega_{\flat}
\colon
\Lambda'_A [sL'']
\longrightarrow
s^2L'
$$
and its composite with the injection
$s^2L' \to T_A [s^2L']$.
\endproclaim

In other words,
we can view 
the morphism
$\iota \Omega_{\sharp}$
as the element
$$
\sum _{i \geq 0} (\Omega_{\flat})^{\cup i} \in
\roman{Hom}_A(\Lambda'_A [sL''],T_A[s^2L'])
\tag3.15.2
$$
of the  graded $A$-algebra
$
\roman{Hom}_A(\Lambda'_A [sL''],T_A[s^2L'])
=
\roman{Alt}_A(L'',T_A[s^2L']) .
$

\demo{Proof of {\rm 3.15}}
This is an immediate consequence of the description
of the induced morphism
of coalgebras
given in (1.6.4.3). \qed
\enddemo

\proclaim{Corollary 3.16}
In the graded $A$-algebra
$\roman{Hom}_A(\Lambda'_A [sL''],T_A[s^2L'])$
we have
$$
D^{\omega}(\iota \Omega_{\sharp}) = 0.
$$
\endproclaim

\demo{Proof} From the generalized Bianchi identity (2.8.2) we know that
$D^{\omega}(\Omega_{\flat}) = 0$.
However, 
$D^{\omega}$
is a derivation over $R$
of the graded $A$-algebra
$\roman{Hom}_A(\Lambda'_A [sL''],T_A[s^2L'])$
since, by construction,
it comes from a derivation
in
$\roman{Alt}_A(L,T_A[s^2L'])$.
Hence we have, for $i \geq 2$, 
$$
D^{\omega}((\Omega_{\flat})^{\cup i}) =
\sum_{j+k = i-1}
(\Omega_{\flat})^{\cup j}
\cup (D^{\omega}(\Omega_{\flat}))
\cup (\Omega_{\flat})^{\cup k}
=0.
$$
Hence
$$
D^{\omega}(\iota \Omega_{\sharp}) = 0. \qed
$$
\enddemo

\demo{Proof of Lemma {\rm 3.11}}
This follows
at once from the fact that
the injection (3.12)
is  compatible with the
operations (3.10) and (3.14)
of covariant derivative. \qed
\enddemo

\demo{Proof of Theorem {\rm 3.8}}
It is clear that the 
$L$-action induces an $L''$-action on the
invariants 
$(\roman{Hom}_A(\Sigma'_A[s^2L'],A)^{L'}$
and that the full invariants
$\roman{Hom}_A(\Sigma'_A[s^2L'],A)^L$
may be rewritten
$$
\roman{Hom}_A(\Sigma'_A[s^2L'],A)^L
=(\roman{Hom}_A(\Sigma'_A[s^2L'],A)^{L'})^{L''}.
$$
However, since
$L'$ acts trivially on $A$,
by adjointness,
we may as well write
$$
\roman{Hom}_A(\Sigma'_A[s^2L'],A)^{L'}
\cong
\roman{Hom}_A(A\otimes_{L'}\Sigma'_A[s^2L'],A)
$$
where as usual
$A\otimes_{L'}\Sigma'_A[s^2L']$
refers to
$\Sigma'_A[s^2L']$
with the $L'$-action divided out,
and the $L$-action on
$\roman{Hom}_A(\Sigma'_A[s^2L'],A)$
passes to an
$L''$-action on
$\roman{Hom}_A(A\otimes_{L'}\Sigma'_A[s^2L'],A)$.
Furthermore, with respect to this action, we have
$$
\roman{Hom}_A(\Sigma'_A[s^2L'],A)^L
\cong(\roman{Hom}_A(A\otimes_{L'}\Sigma'_A[s^2L'],A))^{L''}.
$$
However,
since on
$A\otimes_{L'}\Sigma'_A[s^2L']$
the $L'$-action 
has been divided out,
the covariant derivative
(3.10)
passes to the 
differential
in
$\roman{Alt}_A(L'',A\otimes_{L'}\Sigma'_A[s^2L'])$
associated to the
$L''$-action on
$A\otimes_{L'}\Sigma'_A[s^2L']$
and,
in view of (3.11.1),
the morphism
$\Omega_{\sharp}$
passes to a cycle
(say)
$$
\Omega_{\sharp\sharp}
\colon
\Lambda'_A [sL'']
\longrightarrow
A\otimes_{L'}\Sigma'_A[s^2L']
$$
in
$\roman{Alt}_A(L'',A\otimes_{L'}\Sigma'_A[s^2L'])$,
that is to say, with respect to the differential
$d = d^1$ 
on $\roman{Alt}_A(L'',A\otimes_{L'}\Sigma'_A[s^2L'])$
given by ($\roman{ 1.3}'$)
(spelled out in (1.6.2)), we have
$$
d(\Omega_{\sharp\sharp}) = 0 
\in  \roman{Alt}_A(L'',A\otimes_{L'}\Sigma'_A[s^2L']).
$$
When we now rewrite the Chern-Weil map
(3.8.2) in the form
$$
\Omega_{\sharp\sharp}^*
\colon
\roman{Hom}_A(A\otimes_{L'}\Sigma'_A[s^2L'],A)^{L''}
\longrightarrow
\roman{Alt}_A(L'',A)
$$
we see that
its image in
$\roman{Alt}_A(L'',A)$
consists indeed of cycles only as asserted. 
In fact, an element
$\varphi \in \roman{Hom}_A(A\otimes_{L'}\Sigma'_A[s^2L'],A)^{L''}$
is just a morphism
$$
\varphi \colon A\otimes_{L'}\Sigma'_A[s^2L'] \longrightarrow A
$$
of $(A,L'')$-modules,
and $\Omega_{\sharp\sharp}^*(\varphi)$
coincides with the image 
$\varphi_*(\Omega_{\sharp\sharp}) \in \roman{Alt}_A(L'',A)$
of the cycle
$\Omega_{\sharp\sharp} \in \roman{Alt}_A(L'',A\otimes_{L'}\Sigma'_A[s^2L'])$
under the induced map
$$
\varphi_*
\colon
\roman{Alt}_A(L'',A\otimes_{L'}\Sigma'_A[s^2L'])
\longrightarrow
\roman{Alt}_A(L'',A),
$$
that is, with the composite
$$
\phi\Omega_{\sharp\sharp}
\colon
\Lambda'_A [sL'']
\longrightarrow
A.
$$
\smallskip
Finally,
a different choice $\omega' \colon L'' \to L$
of 
$\roman {\bold e}$-connection yields an
$\roman {\bold e}$-curvature $\Omega'$ so that the two cycles
$
\Omega_{\sharp\sharp}
$
and
$\Omega'_{\sharp\sharp}$
differ by a boundary.
This proves Theorem 3.8. \qed
\enddemo

\smallskip
\noindent
{\smc Remark 3.17.}
The proof given above shows
in particular
that the {\it true\/}
global invariant 
for $\roman {\bold e}$
is the class
$$
[\Omega_{\sharp\sharp}]
\in \roman H^*\left(\roman{Alt}_A(L'',A\otimes_{L'}\Sigma'_A[s^2L'])\right).
$$
\smallskip
We finally relate our
Chern-Weil construction 
with the classical one:
Let $\xi\colon P \to N$ be a principal bundle
with structure group a compact Lie group $G$,
write $A = C^{\infty}(N)$,
let $L'=\fra g(\xi)$,
viewed as an $(\Bobb R,A)$-Lie algebra with trivial action on $A$,
and let $L= \roman E(\xi)$ and $L''=\roman{Vect}(N)$;
cf. (2.2).
Consider the  extension
{\rm (cf. (2.2.2))}
$$
\roman {\bold e}(\xi)
\colon
0
@>>>
L'
@>>>
L
@>>>
L''
@>>>
0
$$
of  $(\Bobb R,A)$-Lie algebras.

\proclaim{Lemma 3.18}
The algebra of $L$-invariants
$\roman{Hom}_A(\Sigma'_A[s^2L'],A)^L$ is as a
graded commutative algebra over the reals
isomorphic to
the $G$-invariants
$\left(\roman{Hom}_{\Bobb R}(\Sigma'_{\Bobb R}[s^2\fra g],\Bobb R)\right)^G$,
with $G$-action on
$\Sigma'_{\Bobb R}[s^2 \fra g]$ induced by the adjoint representation.
Consequently, when $G$ is connected,
the algebra of $L$-invariants
$\roman{Hom}_A(\Sigma'_A[s^2L'],A)^L$ is as a
graded commutative algebra over the reals
isomorphic to the
$\fra g$-invariants
$\left(\roman{Hom}_{\Bobb R}
(\Sigma'_{\Bobb R}[s^2\fra g],\Bobb R)\right)^{\fra g}$,
with $\fra g$-action on
$\Sigma'_{\Bobb R}[s^2 \fra g]$ induced by the adjoint representation
of $\fra g$ on itself.
\endproclaim

In particular,
the algebra of $L$-invariants
$\roman{Hom}_A(\Sigma'_A[s^2L'],A)^L$ 
depends only on $G$ and the adjoint action
(and not on the specific extension $\roman {\bold e}(\xi)$
of Lie-Rinehart algebras, that is, not explicitly on $L$);
further,
the algebra
$\roman{Hom}_{\Bobb R}(\Sigma'_{\Bobb R}[s^2\fra g],\Bobb R)$
is the algebra
$\Bobb R[c_1,\dots,c_m]$
 of polynomials on a basis 
$\{c_1,\dots,c_m\}$
of
$\roman{Hom}_{\Bobb R}([s^2\fra g],\Bobb R)$,
and hence
the algebra of $L$-invariants
$\roman{Hom}_A(\Sigma'_A[s^2L'],A)^L$ is as a
graded commutative algebra over the reals
isomorphic to
the $G$-invariants
$\Bobb R[c_1,\dots,c_m]^G$ and thence,
when $G$ is connected,
to
the $\fra g$-invariants
$\Bobb R[c_1,\dots,c_m]^{\fra g}$.
Thus we see that
$\roman{Hom}_A(\Sigma'_A[s^2L'],A)^L$
is isomorphic to the algebra of $G$-invariants of the algebra of polynomial
functions on $\fra g$ with respect to an $\Bobb R$-basis,
that is, under the present circumstances,
the source of our Chern-Weil map (3.8.2)
already looks like an appropriate classical object.

\demo{Proof of {\rm (3.18)}}
Write $B = C^{\infty}(P)$, 
$L_N=\Gamma(\tau_N)$,
and $L_B=\Gamma(\tau_P)$.
By construction, in view of what was said about
Atiyah sequences in (2.2),
$L=E(\xi)$ coincides with the invariants  $L_B^G$.
Since $\xi$ is a principal  bundle for
$\roman{ad}(\xi)$,
as $B$-modules,
the induced module $B \otimes_A L'$
is that of sections of the induced bundle $\xi^*(\roman{ad}(\xi))$
and hence isomorphic to $B \otimes_{\Bobb R}\fra g$;
under this isomorphism, the $G$-module structure
on $B \otimes_A L'$ coming from the $G$-action on the first factor
$B$ corresponds to the diagonal action of $G$
on $B \otimes_{\Bobb R}\fra g$.
Furthermore,
$\roman{Hom}_A(\Sigma'_A[s^2L'],B)$
being endowed with the obvious $G$-action induced by the 
$G$-action on $B$,
the obvious map
$$
\roman{Hom}_A(\Sigma'_A[s^2L'],A)
\longrightarrow
\roman{Hom}_A(\Sigma'_A[s^2L'],B)^G
$$
into the $G$-invariants
is an
isomorphism.
On the other hand, we have a chain 
$$
\align
\roman{Hom}_A(\Sigma'_A[s^2L'],B)
&\cong
\roman{Hom}_B(B \otimes _A\Sigma'_A[s^2L'],B)
\\
&\cong
\roman{Hom}_B(\Sigma'_B[B \otimes _A s^2L'],B)
\\
&\cong
\roman{Hom}_B(\Sigma'_B[B \otimes s^2\fra g],B)
\\
&\cong
\roman{Hom}_{\Bobb R}(\Sigma'_{\Bobb R}[s^2\fra g],B)
\endalign
$$
of obvious isomorphisms of graded $\Bobb R$-algebras
and, in view of what has been said before,
the resulting isomorphism
$$
\roman{Hom}_A(\Sigma'_A[s^2L'],B)
\longrightarrow
\roman{Hom}_{\Bobb R}(\Sigma'_{\Bobb R}[s^2\fra g],B)
$$
is one of $G$-modules,
provided we take the $G$-structure
on
$\roman{Hom}_{\Bobb R}(\Sigma'_{\Bobb R}[s^2\fra g],B)$
given by
$$
(y,\phi)\longmapsto
y(\phi) = y \circ\phi\circ \roman{Ad}(y^{-1}),
\quad
y \in G,\ \phi \in \roman{Hom}_{\Bobb R}(\Sigma'_{\Bobb R}[s^2\fra g],B),
$$
where, with an abuse of notation, the symbol Ad
refers to the $G$-action on
$\Sigma'_{\Bobb R}[s^2\fra g]$ 
induced by the adjoint representation of $G$ on $\fra g$.
Furthermore, the isomorphic objects
$$
\roman{Hom}_B(B \otimes _A\Sigma'_A[s^2L'],B),
\quad
\roman{Hom}_B(\Sigma'_B[B \otimes _A s^2L'],B),
\quad
\roman{Hom}_B(\Sigma'_B[B \otimes s^2\fra g],B)
$$
admit structures of $(B,L_B)$-modules and of  $G$-modules
in such a way that, if we write $M$ for any of these,
we have
$$
(M^G)^L = (M^{L_B})^G.
\tag3.19
$$
In view of (3.19), we conclude that
$$
\align
\roman{Hom}_A(\Sigma'_A[s^2L'],A)^L
&\cong
\left(\roman{Hom}_B(\Sigma'_B[B \otimes s^2\fra g],B)^G\right)^L
\\
&\cong
\left(\roman{Hom}_B(\Sigma'_B[B \otimes s^2\fra g],B)^{L_B}\right)^G
\\
&\cong
\left(\roman{Hom}_{\Bobb R}(\Sigma'_{\Bobb R}[s^2\fra g],B^{L_B})\right)^G
\\
&\cong
\left(\roman{Hom}_{\Bobb R}(\Sigma'_{\Bobb R}[s^2\fra g],\Bobb R)\right)^G. 
\qed
\endalign
$$
\enddemo

\smallskip
To explain the formal nature of our  argument 
relating our Chern-Weil map 
(3.8.1)
with the classical one,
we momentarily
return to an arbitrary ground ring
$R$, arbitrary commutative $R$-algebra $A$, 
and arbitrary extension (2.1)
of Lie-Rinehart algebras,
subject only to the condition that
the extension split in the category of $A$-modules.
We then consider 
the
graded $A$-algebra
$\Sigma_A^* [s^2L']= \roman{Hom}_A(\Sigma_A[s^2L'],A)$ 
of $A$-multilinear symmetric functions on $s^2L'$ with values in $A$,
with the shuffle product
as multiplication;
cf. what is said in (1.6.4) above.
For $\alpha \in L$ and
$\phi \in 
\roman{Hom}_A(\Sigma_A[s^2L'],A)$,
let
$$
(\zeta(\alpha))(\phi) = \phi \circ (\roman{ad}(\alpha)) - \alpha \circ \phi
\colon
\Sigma_A[s^2L']
\longrightarrow
A;
$$
here $(\alpha \circ \phi)(x) =
\alpha(\phi(x))$
where, with an abuse of notation,
the operator on $A$ induced by $\alpha$
is denoted by $\alpha$ as well.
Inspection shows that this induces an action 
$$
L \to
\roman{End} 
(\Sigma_A^*[s^2L'])
$$
of $L$
on 
$\Sigma_A^*[s^2L']$.
The 
morphism (1.6.4.13) of graded $A$-algebras
now looks like
$$
\roman{Hom}_A(\Sigma'_A[s^2L'],A)
\longrightarrow
\Sigma^*_A[s^2L']
\tag3.20
$$
and is compatible with the $L$-structures
and natural in the data in the appropriate sense;
when $L'$ is free, this map 
sends a polynomial function in a basis of $L'$
to the corresponding symmetric function on $L'$,
cf. what is said in
(1.6.4) above.
\smallskip
We now suppose that
the ground ring $R$ contains the
rationals so that (3.20) is an isomorphism.
Then the composite of the 
inverse
of (3.20),
restricted to the $L$-invariants,
with the Chern-Weil map
(3.1.1),
yields a morphism
$$
(\Sigma_A^*[s^2L'])^L
\to
\roman{Alt}_A(L'',A)
\tag3.21
$$
of differential graded $A$-algebras
which 
is defined
on the invariant symmetric functions
and coincides formally with the classical Chern-Weil map
\cite{\duponboo,\ \kobanomi}.
Under the circumstances spelled out just before (3.18),
in view of the statement of (3.18), we have
$$
(\Sigma^*_A[s^2L'])^L 
\cong \roman{Inv}(\fra g) =
(\Sigma^*_{\Bobb R}[s^2\fra g])^{G} 
$$
and, keeping in mind that
$$
\roman H^*(\roman{Alt}_A(L'',A))
\cong \roman H^*_{\roman{de Rham}}(N,\Bobb R),
$$
our Chern-Weil map (3.8.1) 
boils indeed down to the classical one
$$
\roman{Inv}(\fra g) \longrightarrow \roman H^*_{\roman{de Rham}}(N,\Bobb R),
$$
perhaps up to certain multiplicative factors
(involving the factorials)
which are due to the different possibilities
of defining the de Rham algebra in characteristic zero
and to the definition of the Chern-Weil map directly in terms of
the algebra of symmetric functions.
We note that
the description of the
Chern-Weil map in \cite\milnstas\ 
involves actually
the algebra of polynomial functions
(on $\fra g$)
and not that of symmetric functions
and is considerably closer
to ours than that in the other sources
\cite{\duponboo,\ \kobanomi}.

\beginsection 4. Examples

(1) We take the ground ring
$R$ to be that of the reals $\Bobb R$.
Let $N$ be a smooth finite dimensional manifold, 
write $A = C^{\infty}(N)$,
let
$L'' = \roman{Vect}(N)$, the 
$(\Bobb R,A)$-Lie algebra
of smooth vector fields on $N$,
and let $L'$ be the $(\Bobb R,A)$-Lie algebra
which as an $A$-module is just $A$, with trivial
(= abelian) Lie algebra structure.
Then, on the one hand, 
$\roman H^2(\roman{Alt}_A(L'',A)$
is isomorphic to
$\roman H^2_{\roman{de Rham}}(N,\Bobb R)$,
see 
\cite\rinehart\ 
for details
while, on the other hand,
in view of Theorem 2.6,
the cohomology group
$\roman H^2(\roman{Alt}_A(L'',A)$
classifies extensions of $(\Bobb R,A)$-Lie algebras of the kind
$$
{\roman {\bold e}}
\colon
0 
@>>>
A
@>>>
L
@>>>
L''
@>>>
0
\tag4.1
$$
subject to the requirement that
the adjoint representation of $L$ on itself induces the
$L''$-module structure on $A$.
Under these circumstances our Chern-Weil map (3.8.2)
has a very simple form.  Indeed,
as an $A$-module, $L'$ may be written as an induced module
$L' = A \otimes_{\Bobb R} \Bobb R$, and 
the same is manifestly true of the graded symmetric coalgebra
$\Sigma'_A[s^2L']$ which is in fact isomorphic to
$A \otimes _{\Bobb R}\Sigma'_{\Bobb R}[s^2{\Bobb R}]$.
Hence the graded commutative algebra
$\roman {Hom}_A(\Sigma'_A[s^2L'],A)$
may be rewritten
$\roman {Hom}_{\Bobb R}(\Sigma'_{\Bobb R}[s^2\Bobb R],A)$,
and hence
the subalgebra of invariants 
$\roman {Hom}_A(\Sigma'_A[s^2L'],A)^L$
looks like
$\roman {Hom}_{\Bobb R}(\Sigma'_{\Bobb R}[s^2\Bobb R],\Bobb R)$,
which is just the polynomial algebra
$\Bobb R[c]$ on a basis $\{c\}$ of
$\roman {Hom}_{\Bobb R}(s^2\Bobb R,\Bobb R)$.
Consequently our Chern-Weil map
(3.8.2) looks like
$$
\Bobb R[c]
\longrightarrow
\roman H^{2*}_{\roman{de Rham}}(N,\Bobb R).
\tag4.2
$$
When the extension
${\roman {\bold e}}$
does not come from a principal $S^1$-bundle,
that is to say, when the class
$[\Omega]\in \roman H^2_{\roman{de Rham}}(N,\Bobb R)$
is not integral
(i.~e. does not have integral periods),
this kind of example is not covered
by the classical theory.
\smallskip\noindent
(2) More generally,
let
$\xi \colon P \to N$ 
be a principal bundle,
with structure group $G$ and Lie algebra $\fra g$, 
let
$L' = \fra g(\xi) = \Gamma(\roman{ad}(\xi))$ be the 
$A$-Lie algebra, viewed as an
$(\Bobb R,A)$-Lie algebra
with trivial action on $A$,
which as an $A$-module
is the space of sections of the adjoint 
bundle $\roman{ad}(\xi)$, and consider an 
extension
$$
\roman {\bold e}
\colon
0 
@>>>
L'
@>>>
L
@>>>
L''
@>>>
0
\tag4.3
$$
of Lie-Rinehart algebras
having the same outer action
$L'' \to \roman{ODer}(L')$ as that coming from  the
extension
$$
\roman {\bold e(\xi)}
\colon
0 
@>>>
L'
@>>>
E(\xi)
@>>>
L''
@>>>
0
$$
given in (2.2.2); we remind the reader that the notion of outer action
has been reproduced in Section 2.
Under these circumstances,
we can  still conclude that
the subalgebra of invariants 
$\roman {Hom}_A(\Sigma'_A[s^2L'],A)^L$
looks like
$\roman {Hom}_{\Bobb R}(\Sigma'_{\Bobb R}[s^2 \fra g],\Bobb R)^G$,
that is to say, it is the algebra of $G$-invariants
of the algebra 
$\Bobb R[c_1,c_2,\dots,c_m]$
of polynomial functions on $\fra g$ with respect to an
$\Bobb R$-basis
$\{c_1,c_2,\dots,c_m\}$
of
$\roman {Hom}_{\Bobb R}(s^2 \fra g,\Bobb R)$.
In fact,
the 
$L$-action on
$\roman {Hom}_A(\Sigma'_A[s^2L'],A)$
passes to an $L''$-action on
the
algebra  of invariants
$\roman {Hom}_A(\Sigma'_A[s^2L'],A)^{L'}$,
determined entirely by the
corresponding outer action
$L'' \to \roman{ODer}(L')$,
and the
algebra  of invariants
$\roman {Hom}_A(\Sigma'_A[s^2L'],A)^L$
may be rewritten
$$
\roman {Hom}_A(\Sigma'_A[s^2L'],A)^L
=
\left(\roman {Hom}_A(\Sigma'_A[s^2L'],A)^{L'}\right)^{L''}.
$$
However, 
the
algebra  of invariants
$\left(\roman {Hom}_A(\Sigma'_A[s^2L'],A)^{L'}\right)^{L''}$,
in turn, may be rewritten
$$
\left(\roman {Hom}_A(\Sigma'_A[s^2L'],A)^{L'}\right)^{L''}
=
\roman {Hom}_A(\Sigma'_A[s^2L'],A)^{E(\xi)},
$$
and, by virtue of (3.18), we know that,
as a graded commutative algebra over the reals,
the algebra $\roman {Hom}_A(\Sigma'_A[s^2L'],A)^{E(\xi)}$
is 
isomorphic to the algebra of invariants
$\roman {Hom}_{\Bobb R}(\Sigma'_{\Bobb R}[s^2 \fra g],\Bobb R)^G$.
Consequently our Chern-Weil map
(3.8.2) looks like
$$
\Bobb R[c_1,c_2,\dots,c_m]^G
\longrightarrow
\roman H^{2*}_{\roman{de Rham}}(N,\Bobb R).
\tag4.4
$$
When the extension
$\roman {\bold e}$
does not come from a principal $G$-bundle,
again this kind of example is not covered
by the classical theory.
To obtain explicit examples,
suppose that the Lie algebra $\fra g$ has a non-trivial
centre $\fra z$,
and let
$\zeta (\xi) \colon P \times _G \fra z \to N$
be the corresponding associated bundle with fibre $\fra z$---
this bundle is trivial when $G$ is connected.
When
the cohomology group
$\roman H^2_{\roman{de Rham}}(N,\zeta(\xi))$
is non-zero,
Theorem 2.7
provides a wealth of examples
of extensions
of 
$L''$ by $L'$
which do not come from a principal bundle
but have
the same outer action
of $L''$ on $L'$
as that coming from the 
principal bundle $\xi$ we started with.
An approach to the corresponding global theory, phrased in terms of
Lie groupoids, has been given 
in \cite\mackthr.
\smallskip\noindent
(3) Let $N$ be a smooth finite dimensional manifold,
let $C \subseteq N$ be a compact subset, not necessarily a smooth manifold,
and let $A_C$ be the algebra of smooth functions on $C$ in the sense of
{\smc Whitney\/}~\cite\whitnone,~\cite\whitntwo.
It will here be convenient to take
for a {\it smooth function\/}
$f$ on $C$ in this sense  a class of smooth functions
$h$ defined on $N$, 
two functions being identified
whenever they coincide on $C$.
(It is also customary to take
classes of smooth functions
$h$ defined only on a neighborhood of $C$ in $N$.)
Let $I_C$ be the ideal of smooth functions on $N$
that vanish on $C$, so that
$A_C= C^{\infty}(N)/I_C$.
Furthermore, let $\roman{Vect}(N,C) \subseteq \roman{Vect}(N)$
be the set of smooth vector fields $X$ on $N$ that preserve $I_C$
in the sense that
$$
Xh = 0 \quad\text{on $C$ whenever}\quad h = 0 \quad\text{on}\quad C .
$$
It is readily seen that
$\roman{Vect}(N,C)$ inherits a structure of an
$(\Bobb R, C^{\infty}(N))$-Lie algebra
from
$\roman{Vect}(N)$.
Let $L_C =  A_C\otimes_{C^{\infty}(N)} \roman{Vect}(N,C)$;
inspection shows that the 
$(\Bobb R, C^{\infty}(N))$-Lie algebra structure on
$\roman{Vect}(N,C)$
passes to that of an
$(\Bobb R, A_C)$-Lie algebra
on $L_C$.
In the special case where $C$ is a smooth submanifold of $N$,
the obvious  map
$$
\roman{Vect}(N,C)
\longrightarrow
\roman{Vect}(C)
$$
induces an isomorphism
$$
L_C \longrightarrow
\roman{Vect}(C).
$$
Hence we refer to $L_C$ as the
$(\Bobb R,A_C)$-Lie algebra of {\it smooth vector fields
on\/} $C$. 
\smallskip
As in (2),
write $A=C^{\infty}(N)$,
let 
$\xi \colon P \to N$ 
be a principal bundle,
with structure group $G$ and Lie algebra $g$, 
and consider the
$A$-Lie algebra
$\fra g(\xi)$, viewed as an
$(\Bobb R,A)$-Lie algebra
with trivial action on $A$.
Let $\fra g_C(\xi)=  A_C\otimes_{C^{\infty}(N)}  \fra g(\xi)$;
it 
is manifestly a projective $A_C$-module and
inherits an obvious structure of an
$A_C$-Lie algebra;
it will henceforth be viewed as an
$(\Bobb R,A_C)$-Lie algebra
with trivial action on $A_C$.
Under these circumstances,
our Chern-Weil map
(3.8.2) arising from an arbitrary
extension
$$
\roman {\bold e}
\colon
0 
@>>>
\fra g_C(\xi)
@>>>
L
@>>>
L_C
@>>>
0
\tag4.5
$$
of Lie-Rinehart algebras
looks like
$$
\roman {Hom}_{A_C}(\Sigma'_{A_C}[s^2 \fra g_C(\xi)],A_C)^{L}
\longrightarrow
\roman H^{2*}_{A_C}\left(\roman{Alt}_{A_C}(L_C,A_C)\right).
\tag4.6
$$
Again this kind of examples is not covered by the
classical approach.
In particular, 
we may restrict the corresponding sequence (2.2.2) 
to
the $(\Bobb R,C^{\infty}(N))$-Lie algebra
$\roman{Vect}(N,C)$, as indicated in the commutative diagram
$$
\CD
0
@>>>
\fra g(\xi)
@>>>
\roman E(\xi,C)
@>>>
\roman{Vect}(N,C)
@>>>
0
\\
@.
@VVV
@VVV
@VVV
@.
\\
0
@>>>
\fra g(\xi)
@>>>
\roman E(\xi)
@>>>
\roman{Vect}(N)
@>>>
0
\endCD
\tag4.7
$$
where
the
$(\Bobb R,C^{\infty}(N))$-Lie algebra
$\roman E(\xi,C)$ is 
defined by the requirement that
$\roman E(\xi,C)$,
$\roman E(\xi)$,
$\roman{Vect}(N,C)$, and
$\roman{Vect}(N)$
constitute a pull back diagram.
Since, as
$(C^{\infty}(N))$-modules,
the bottom row of (4.8) splits, so does the  top row;
consequently, with the notation
$\roman E_C(\xi) = A_C\otimes_{C^{\infty}(N)}\roman E(\xi,C)$,
the corresponding sequence
$$
\CD
0
@>>>
\fra g_C(\xi)
@>>>
\roman E_C(\xi)
@>>>
L_C
@>>>
0
\endCD
\tag4.8
$$
is still exact and in particular an extension
of
$(\Bobb R,A_C)$-Lie algebras.
Our Chern-Weil map (4.6)
furnishes characteristic classes for it.
\smallskip
More generally,
we now 
consider a general extension
of Lie-Rinehart algebras
of the kind (4.5), subject only to the condition that
its corresponding  outer action
$$
L_C \longrightarrow \roman{ODer}(\fra g_C(\xi)),
$$
cf. Section 2, coincides with that
for the extension (4.8).
The argument already used before shows that
the algebra  of invariants
$\roman {Hom}_{A_C}(\Sigma'_{A_C}[s^2 \fra g_C(\xi)],A_C)^{L}$
may be rewritten
$$
\align
\roman {Hom}_{A_C}(\Sigma'_{A_C}[s^2 \fra g_C(\xi)],A_C)^{L}
&=
\left(\roman {Hom}_{A_C}
(\Sigma'_{A_C}[s^2 \fra g_C(\xi)],A_C)^{L'}\right)^{L_C},
\\
&=\roman {Hom}_{A_C}(\Sigma'_{A_C}[s^2 \fra g_C(\xi)],A_C)^{E_C(\xi)},
\endalign
$$
where $L'= \fra g_C(\xi)$,
and hence does not depend on the extension (4.5) but only
on the corresponding outer action.
Furthermore,
restriction induces a commutative diagram
$$
\CD
\Bobb R[c_1,c_2,\dots,c_m]^G
@>>>
\roman H^{2*}_{\roman{de Rham}}(N,\Bobb R)
\\
@VVV
@VVV
\\
\roman {Hom}_{A_C}(\Sigma'_{A_C}[s^2 \fra g_C(\xi)],A_C)^{E_C(\xi)}
@>>>
\roman H^{2*}(\roman{Alt}_{A_C}(L_C,A_C)),
\endCD
\tag4.9
$$
the horizontal arrows being the corresponding
Chern-Weil maps.
\smallskip
This discussion raises, among others, the following two questions:
\roster
\item"(i)"
Is the vertical morphism
$$
\Bobb R[c_1,c_2,\dots,c_m]^G
\longrightarrow
\roman {Hom}_{A_C}(\Sigma'_{A_C}[s^2 \fra g_C(\xi)],A_C)^{E_C(\xi)}
$$
in (4.9) an isomorphism?
\item"(ii)"
Can we intrinsically define the notion of a smooth principal bundle
for $A_C$
merely over $C$ which is {\it not\/}
necessarily induced from a principal bundle over $N$?
\item"(iii)"
If the answer to (2) is yes,
can we then determine 
the corresponding algebra of invariants
directly in terms of the new structure over $C$?
\endroster
\smallskip\noindent
(4) Let $\Cal F$ be a {\it transversally complete\/}
foliation of a smooth manifold $V$ \cite\almolone,
write $\tau_{\Cal F} \colon T\Cal F \to V$
for the tangent bundle of $\Cal F$,
and let $\nu_{\Cal F} \colon Q \to V$
be its normal bundle, so that $Q = TV /T\Cal F$.
Let $E(\Cal F)$ be the Lie algebra of vector fields
on $V$ preserving the foliation,
and let $A$ be the algebra of smooth functions
on the leaf space, i.~e. smooth functions
on $V$ which are
constant on the leaves.
Then, with the obvious structure, the pair $(A,E(\Cal F))$
constitutes a Lie-Rinehart algebra.
Since
$\Cal F$
is transversally complete,
the closures of the leaves constitute a smooth fibre bundle
$F \to V \to W$,
and the algebra $A$ may be identified with the algebra of smooth
functions on $W$;
moreover, 
the obvious map
from $E(\Cal F)$
to $\roman{Vect}(W)$
which is part of the Lie-Rinehart structure
of
$(A,E(\Cal F))$
is surjective,
and there results
an extension
$$
\roman {\bold e_{\Cal F}}
\colon
0 
@>>>
L'
@>>>
E(\Cal F)
@>>>
\roman{Vect}(W)
@>>>
0
\tag4.10
$$
of Lie-Rinehart algebras.
The kernel
$L'$ is in fact the space of sections of a Lie algebra
bundle on $W$.
Our Chern-Weil construction
yields characteristic classes 
in 
$\roman H^{2*}_{\roman{de Rham}}(W,\Bobb R)$
for this extension.
\smallskip
We conclude with an illustration
which I learnt from A. Weinstein:
Let $V = \roman{SU}(2) \times \roman{SU}(2)$,
and let $\Cal F$ be the foliation
defined by a dense one-parameter subgroup
in the maximal torus $S^1 \times S^1$
in
$\roman{SU}(2) \times \roman{SU}(2)$.
Then the space
$W$ is $S^2 \times S^2$,
and the Chern-Weil construction yields a characteristic class
in
$\roman H^2_{\roman{de Rham}}(S^2 \times S^2,\Bobb R)$
which may be viewed as an {\it irrational Chern  class\/}.
In view of a result
of {\smc Almeida and Molino} \cite\almolone,
the transitive Lie algebroid
corresponding to
(4.10)
does not integrate to a principal bundle;
in fact, {\smc Mackenzie's} integrability obstruction \cite\mackbook\ 
is non-zero.
It is clear that there are many other examples of this kind.

\medskip
\centerline{References}
\smallskip
\
\widestnumber\key{999}
\ref \no \almolone
\by R. Almeida and P. Molino
\paper Suites d'Atiyah et feuilletages transversalement compl\`etes
\jour C. R. Acad. Sci. Paris I 
\vol 300
\yr 1985
\pages 13--15
\endref
\ref \no \atiyaone
\by M. F. Atiyah
\paper Complex analytic connections in fibre bundles
\jour Trans. Amer. Math. Soc.
\vol 85
\yr 1957
\pages 181--207
\endref
\ref \no \cartanon
\by H. Cartan
\paper Notions d'alg\`ebre diff\'erentielle;
applications aux groupes de Lie et aux vari\'et\'es
o\`u op\`ere un groupe
de Lie
\jour Coll. Topologie Alg\'ebrique
\vol 
\paperinfo Bruxelles
\yr 1950
\pages  15--28
\endref
\ref \no \cartanse
\by  H. Cartan
\paper Alg\`ebres d'Eilenberg--Mac Lane et homotopie
\paperinfo expos\'e 7: Puissances divis\'ees
\jour S\'eminaire H. Cartan 1954/55
\publ Ecole Normale Superieure, Paris, 1956
\endref
\ref \no \carteile
\by H. Cartan and S. Eilenberg
\book Homological Algebra
\publ Princeton University Press
\publaddr Princeton
\yr 1956
\endref
\ref \no \cheveile
\by C. Chevalley and S. Eilenberg
\paper Cohomology theory of Lie groups and Lie algebras
\jour  Trans. Amer. Math. Soc.
\vol 63
\yr 1948
\pages 85--124
\endref
\ref \no \duponboo
\by J. L. Dupont
\book Curvature and characteristic classes
\bookinfo Lecture Notes in Mathematics,
 No. 640
\publ Springer
\publaddr Berlin--G\"ottingen--Heidelberg
\yr 1978
\endref
\ref \no \eilmacon
\by S. Eilenberg and S. Mac Lane
\paper Cohomology theory in abstract groups. I.
\jour Ann. of Math.
\vol 48
\yr 1947
\pages  51--78
\moreref
\paper Cohomology theory in abstract groups. II. 
Group extensions with a non-abelian kernel
\jour Ann. of Math.
\vol 48
\yr 1947
\pages  326--341
\endref

\ref \no \grehalva
\by W. Greub, S. Halperin, and R. Vanstone
\book Connections, Curvature, and Cohomology, vol.s I. -- III.
\bookinfo Pure and Applied Mathematics,
a series of monographs and textbooks
\publ Academic Press
\publaddr New York - San Francisco - London
\yr 1976
\endref

\ref \no \herzone
\by J. Herz
\paper Pseudo-alg\`ebres de Lie
\jour C. R. Acad. Sci. Paris 
\vol 236
\yr 1953
\pages 1935--1937
\endref

\ref \no \poiscoho
\by J. Huebschmann
\paper Poisson cohomology and quantization
\jour J. f\"ur die Reine und Angew. Math.
\vol 408
\yr 1990
\pages 57--113
\endref
\ref \no  \souriau
\by J. Huebschmann
\paper On the quantization of Poisson algebras
\paperinfo Symplectic Geometry and Mathematical Physics,
Actes du colloque en l'honneur de Jean-Marie Souriau,
P. Donato, C. Duval, J. Elhadad, G.M. Tuynman, eds.;
Progress in Mathematics, Vol. 99
\publ Birkh\"auser
\publaddr Boston $\cdot$ Basel $\cdot$ Berlin
\yr 1991
\pages 204--233
\endref
\ref \no \exteweil
\by J. Huebschmann
\paper 
Higher homotopies and perturbations
for extensions of Lie-\linebreak
Rinehart algebras
\paperinfo in preparation
\endref
\ref \no \perturba
\by J. Huebschmann
\paper Perturbation theory and free resolutions for nilpotent
groups of class 2
\jour J. of Algebra
\yr 1989
\vol 126
\pages 348--399
\endref
\ref \no \husmosta
\by D. Husemoller, J. C. Moore, and J. D. Stasheff
\paper Differential homological algebra and homogeneous spaces
\jour J. of Pure and Applied Algebra
\vol 5
\yr 1974
\pages  113--185
\endref

\ref \no \kamtonth
\by F. W. Kamber and Ph. Tondeur
\paper Weil algebras and characteristic classes
for foliated bundles in Cech cohomology
\jour Diff. Geom. Proc. Sym. Pure Math.
\vol XXVII
\yr 1975
\pages 283--294 
\endref

\ref \no \kobanomi
\by S. Kobayashi and K. Nomizu
\book Foundations of differential geometry, I (1963), II (1969)
\bookinfo Interscience Tracts in Pure and Applied Mathematics, No. 15
\publ Interscience Publ.
\publaddr New York-London-Sydney
\endref
\ref \no \koszulon
\by J. L. Koszul
\paper Crochet de Schouten-Nijenhuis et cohomologie
\jour Ast\'erisque,
\vol hors-s\'erie,
\yr 1985
\pages 251--271
\paperinfo in E. Cartan et les Math\'ematiciens d'aujourd'hui, 
Lyon, 25--29 Juin, 1984
\endref
\ref \no \koszultw
\by J. L. Koszul
\book Lectures on fibre bundles and differential geometry
\publ Tata Institute of Fundamental research
\publaddr Bombay
\yr 1960
\finalinfo reprinted: Springer 1986
\endref
\ref \no \maclaboo
\by S. Mac Lane
\book Homology
\bookinfo Die Grundlehren der mathematischen Wissenschaften
 No. 114
\publ Springer
\publaddr Berlin--G\"ottingen--Heidelberg
\yr 1963
\endref
\ref \no \mackbook
\by K. Mackenzie
\book Lie groupoids and Lie algebroids in differential geometry
\bookinfo London Math. Soc. Lecture Note Series, vol. 124
\publ Cambridge University Press
\publaddr Cambridge, England
\yr 1987
\endref
\ref \no \macktwo
\by K. Mackenzie
\paper On extensions of principal bundles
\jour Annals of Global Analysis and Geometry
\vol 6
\yr 1988
\pages 141--163
\endref
\ref \no \mackthr
\by K. Mackenzie
\paper Classifications of principal bundles and Lie groupoids
with prescribed gauge group bundle
\jour J. of Pure and Applied Algebra
\vol 58
\yr 1989
\pages 181--208
\endref
\ref \no \maybook
\by J. P. May
\book Simplicial objects in algebraic topology
\publ van Nostrand
\publaddr Princeton
\yr 1967
\endref
\ref \no \milnstas
\by J. Milnor and J. D. Stasheff
\book Characteristic classes
\bookinfo Annals of Mathematics Studies, no. 76
\publ Princeton University Press
\publaddr Princeton, New Jersey
\yr 1974
\endref
\ref \no \munkholm
\by H. J. Munkholm
\paper The Eilenberg--Moore spectral sequence and strongly homotopy
multiplicative maps
\jour J. of Pure and Applied Algebra
\vol 9
\yr 1976
\pages  1--50
\endref

\ref \no \osborn
\by H. Osborn
\paper The Chern-Weil construction
\jour Diff. Geom. Proc. Symp. Pure Math.
\vol XXVII
\yr 1975
\pages 383--395 
\endref

\ref \no \palaione
\by R. Palais
\paper The cohomology of Lie rings
\jour  Proc. Symp. Pure Math.
\vol III
\yr 1961
\pages 130--137
\paperinfo Amer. Math. Soc., Providence, R. I.
\endref

\ref \no \pradione
\by J. Pradines
\paper Th\'eorie de Lie pour les groupo\"\i des diff\'erentiables.
Relations entre propri\'et\'es locales et globales
\jour C. R. Acad. Sci. Paris 
\vol 263
\yr 1966
\pages 907--910
\endref
\ref \no \praditwo
\by J. Pradines
\paper 
Th\'eorie de Lie pour les groupo\"\i des diff\'erentiables.
Calcul diff\'erentiel dans la cat\'egorie des
groupo\"\i des  infinit\'esimaux
\jour C. R. Acad. Sci. Paris 
\vol 264
\yr 1967
\pages 245--248
\endref
\ref \no \pradithr
\by J. Pradines
\paper 
G\'eom\'etrie diff\'erentielle 
au-dessus d'un groupo\"\i de
\jour C. R. Acad. Sci. Paris 
\vol 266
\yr 1968
\pages 1194--1196
\endref
\ref \no \pradifou
\by J. Pradines
\paper 
Troisi\`eme th\'eor\`eme de Lie
pour les
groupo\"\i des diff\'erentiables 
\jour C. R. Acad. Sci. Paris 
\vol 267
\yr 1968
\pages 21--23
\endref

\ref \no \rinehart
\by G. Rinehart
\paper Differential forms for general commutative algebras
\jour  Trans. Amer. Math. Soc.
\vol 108
\yr 1963
\pages 195--222
\endref

\ref \no \stashtwo
\by J. D. Stasheff
\paper Constrained Hamiltonians
\book Conference on elliptic curves and modular forms in algebraic
topology
\bookinfo Institute for Advanced Study,
Princeton, September 15--17, 1986,
ed. P. S. Landweber,
Lecture Notes in Mathematics, No. 1326
\publ Springer
\publaddr Berlin-Heidelberg- 
New York
\yr 1988
\pages 150--160
\endref

\ref \no \stashfiv
\by J. D. Stasheff
\paper Homological reduction of constrained Poisson algebras
\jour J. of Diff. Geom. (to appear)
\endref

\ref \no \stashnin
\by J. D. Stasheff
\paper Deformation theory and the Batalin-Vilkovisky 
master equation
\paperinfo in: Ascona, 1996 (to appear)
\endref

\ref \no \teleman
\by N. Teleman
\paper A characteristic ring of a Lie algebra extension
\jour Accad. Naz. Lincei. Rend. Cl. Sci. Fis. Mat. Natur. (8)
\vol 52
\yr 1972
\pages 498--506 and 708--711
\endref
\ref \no \whitnone
\by H. Whitney
\paper Analytic extensions of differentiable functions defined
on closed sets
\jour Trans. Amer. Math. Soc.
\vol 36
\yr 1934
\pages  63--89
\endref
\ref \no \whitntwo
\by H. Whitney
\paper On ideals  of differentiable functions 
\jour  Amer. J. of Math. 
\vol 70
\yr 1948
\pages  635--658
\endref
\enddocument